\documentclass[5p]{elsarticle}

 \PassOptionsToPackage{force}{filehook} 

\usepackage[table]{xcolor}
\usepackage[utf8]{inputenc}
\usepackage[T1]{fontenc}
\usepackage{soul}

\usepackage{geometry}
\usepackage{natbib}

\usepackage{cuted}
\usepackage{amsmath}
\usepackage[capitalise]{cleveref}
\usepackage{mathtools}
\usepackage{amsfonts}
\usepackage{multicol}
\usepackage{graphicx}
\usepackage{bbm}
\usepackage{etoolbox}
\usepackage{eurosym}
\usepackage{amssymb}
\usepackage{IEEEtrantools}
\usepackage{breqn}

\usepackage{pgf}
\usepackage{pgfplots}
\usepackage{standalone}
\usepackage{float}
\usepackage{dblfloatfix}    

\IEEEeqnarraydefcolsep{0}{\leftmargini} 
\usepackage[T3,T1]{fontenc}
\DeclareSymbolFont{tipa}{T3}{cmr}{m}{n}
\DeclareMathAccent{\invbreve}{\mathalpha}{tipa}{16}
\usepackage{stmaryrd}
\usepackage{ulem}
\pgfplotsset{compat=newest}
\usetikzlibrary{plotmarks}
\usetikzlibrary{arrows.meta}
\usepgfplotslibrary{patchplots}
\usepgfplotslibrary{units}
\usetikzlibrary{patterns}
\usetikzlibrary{calc,positioning}

\usepackage{booktabs}
\usepackage{multirow}
\usepackage{url}

\newlength\fheight 
\newlength\fwidth

\graphicspath{ {Figures/} {figure/} {Figs/} {Matlab/}}

\makeatletter
\def\input@path{{./Figures/}{Figs/}{Matlab/}}
\makeatother

\usepackage{framed} 

\usepackage[noprefix]{nomencl} 

\newcommand{\txtred}[1]{\textcolor{black}{#1}}

\newcommand{\nomunit}[1]{ \hfill{ #1}}

\newcommand\Nomenclature[3][X]{\nomenclature[#1#3]{#2}{#3}}
\newcommand{\RN}[1]{%
	\textup{\uppercase\expandafter{\romannumeral#1}}%
}
\renewcommand\nomgroup[1]{%
	\item[\bfseries
	\ifstrequal{#1}{A}{Acronyms:}{%
		\ifstrequal{#1}{Q}{Indices/Sets:}{%
			\ifstrequal{#1}{S}{Subsets:}{%
				\ifstrequal{#1}{U}{Superscripts:}{%
					\ifstrequal{#1}{P}{Parameters:}{%
						\ifstrequal{#1}{C}{Continious Variables:}{%
							\ifstrequal{#1}{D}{Discerte Variables:}{}}}}}}}]%
			}				

\makenomenclature
\renewcommand*\nompreamble{\begin{multicols}{2}} 
	\renewcommand*\nompostamble{\end{multicols}}

\usepackage[english]{babel}

  \date{}
  
\begin{document}

\sloppy

  \begin{frontmatter}
  \title{Demand-side management via optimal production scheduling in power-intensive industries: The case of metal casting process  }

\author[rvt]{D.~Ramin}
\ead{danial.ramin@itia.cnr.it}
\author[rvt]{S.~Spinelli \corref{cor1}}
\ead{stefano.spinelli@itia.cnr.it}
\author[rvt]{A.~Brusaferri }
\ead{alessandro.brusaferri@itia.cnr.it}
\fntext[]{Ramin, D., S. Spinelli, and A. Brusaferri. "Demand-side management via optimal production scheduling in power-intensive industries: The case of metal casting process." Applied Energy 225 (2018): 622-636. DOI: 10.1016/j.apenergy.2018.03.084	\textregistered 2018. This manuscript version is made available under the CC-BY-NC-ND 4.0 license \url{http://creativecommons.org/licenses/by-nc-nd/4.0/}}

\cortext[cor1]{Corresponding author }
\address[rvt]{ITIA-CNR, Institute of Industrial Technology and Automation - Consiglio Nazionale delle Ricerche, via Corti 12, 20133 Milano, Italy}

\begin{abstract}
The increasing challenges to the grid stability posed by the penetration of renewable energy resources urge a more active role for demand response programs as viable alternatives to a further expansion of peak power generators. This work presents a methodology to exploit the demand flexibility of energy-intensive industries under Demand-Side Management programs in the energy and reserve markets. To this end, we propose a novel scheduling model for a multi-stage multi-line process, which incorporates both the critical manufacturing constraints and the technical requirements imposed by the market. Using mixed integer programming approach, two optimization problems are formulated to sequentially minimize the cost in a day-ahead energy market and maximize the reserve provision when participating in the ancillary market. The effectiveness of day-ahead scheduling model has been verified for the case of a real metal casting plant in the Nordic market, where a significant reduction of energy cost is obtained. Furthermore, the reserve provision is shown to be a potential tool for capitalizing on the reserve market as a secondary revenue stream.
\end{abstract}


\end{frontmatter}
\section{Introduction}

The last decade has witnessed a major paradigm shift in EU energy market and policies.
In 2014 European leaders adopted a climate and energy framework to ensure a 40\% cut in greenhouse emissions from 1990 level by 2030.  
In addition, the framework sets a binding target to increase the share of renewables to 27\% of the final energy consumption at EU level.
One of the main barriers to the integration of renewable energy sources (RES) is their intermittency and unpredictability. Their non-responsive nature makes the already challenging task of maintaining supply-demand balance even more difficult and, if not properly managed, could jeopardize the grid reliability.\\
In Europe, Transmission System Operators (TSOs) are in charge of ensuring the physical balance of power in the grid.
Lacking cheap and efficient storage systems, TSO traditionally relies on dispatchable fossil fuel generation sources to bring production and demand into balance and stabilize the grid frequency.
In the event of an imbalance, grid operator sequentially calls upon three types of generation reserves - categorized based on their response time \cite{ENTSO-E2004} - to bring the grid frequency back to its nominal value of 50 Hz. In order to provide up/downward regulating reserve, involved power plants need to run slightly under/over their max/minimum generation capacity. Consequently, growing imbalance caused by widespread integration of RES requires new solutions beyond traditional dispatchable resources. A viable alternative to peak generation capacity is engaging consumers in power balancing and using their flexibility to avoid the peaks. Techniques involving such practices are referred to as Demand-Response (DR) or Demand-Side Management (DSM) solutions. \\
Electricity market liberalization can be regarded as the most important enabler for DR initiatives. In Europe DR access to the deregulated markets is granted through balance responsible parties (BRP) and aggregators. BRPs are financially responsible to balance the expected electricity generation and demand profile for the suppliers and consumers under their jurisdiction.
They bid their power profile in the wholesale day-ahead market and subsequently on the real-time market to minimize the deviation of real profile from the contracted one.\\ 
While the regulating market is the main platform for restoring grid balance by dispatching reserve capacity, day-ahead market plays a crucial proactive role in the power-balancing decisions taken on a system-wide scale.\\ 
The participation of loads in the day-ahead energy and reserve markets - respectively referred to as price-based and incentive-based programmes in DR literature - creates a win-win situation for both the TSO and end user. While it helps TSO in balancing the power via BRPs, consumers can increase their economic welfare by exploiting their consumption flexibility. 
Furthermore, the consequent growth in the penetration of RES lowers the marginal electricity price by pushing the more expensive fossil fuel power plant out of the market, which is described as "merit-order effect" \cite{Sensfuss2008}.\\
Nevertheless, many flexible loads cannot access these markets due to regulatory and technical barriers, the most prohibitive being costly subscription fees for energy and minimum bid size for the regulating power market \cite{Biegel2014}. This has led to the emergence of aggregators as DR providers for small consumers \cite{shen2014role}. 
Accordingly, in the DR literature, two types of flexible consumers can be identified: small distributed loads such as residential or commercial consumers and single energy-intensive industrial units that are qualified to individually participate in energy or capacity market. 
In the former group, flexibility is mostly provided by shifting or shedding small non-critical loads such as heating and cooling. Even though such loads are straightforward to manage at the user level, there are still many challenges regarding their aggregation such as complexity of distributed control strategies at the aggregator level and the implementation of advanced metering and communication infrastructure and management systems. This has been extensively addressed in recent years under the "smart grid" paradigm \cite{Siano2014}. \\
Industrial loads, on the other hand, can deliver a much larger flexibility without aggregation, even though industrial clusters may also adopt such a paradigm to exploit their collective flexibility \cite{aryandoust2017potential}.
\subsection{Literature Review}
The potential benefits of industrial Demand-Side Management (iDSM) has been noticed both by academia and industry \cite{Paulus2011,Samad2012,Merkert2015}. However, what makes the DSM of industrial consumers challenging is rather the complexity of their underlying processes which demands a deep domain knowledge. \\
In a manufacturing process, in addition to product delivery commitments, there are many other critical constraints regarding the operational feasibility and safety that need to be taken into account when providing load flexibility \cite{wang2013time}. There is a large body of literature dedicated to these problems in the context of scheduling and planning.\\
Short-term scheduling in manufacturing process concerns resource allocation, task timing, and sequencing to attain production objectives.
Various optimization techniques have been reported in the literature that can solve the scheduling problems effectively. A comprehensive review of such techniques can be found in \cite{Harjunkoski2014}.
Among these approaches, mathematical programming, and in particular, Mixed-Integer Linear Programming (MILP) techniques have been widely known to be the most effective ones considering the size, complexity, and flexibility as the key attributes of most manufacturing processes \cite{Merkert2015}. \citet{Velez2014} have summarized the recent advances in MILP techniques for scheduling of large-scale production systems.\\
\txtred{On the modelling front, the most important decision is whether to model the time as a continuous or discrete variable. In discrete-time representation all the events  are confined to take place at finite predefined instances on a time grid. By discretization, part of the solution space is ignored, introducing a  trade-off between solution quality and problem size.
In contrast to continuous-time models, constraints such as material or resource balance need to be monitored only at specific time instances reducing the model complexity, particularly when the underlying process includes inventory and energy resources.	
In addition, discrete-time models are more readily extendible to handle the discrete nature of contractual conditions for trading energy and capacity. \\}

Even though many of the energy-intensive manufacturing processes studied in literature are considered ideal candidates for demand-response programmes, relatively few efforts in scheduling and planning communities have been directed to address market interaction with process scheduling and its modelling aspects. \\
One of the earliest works in this regard is a production planning method proposed by \citet{Daryanian1989} for air separation plant which could take into account spot electricity price. They presented how storage capacity can provide flexibility to minimize electricity cost without service curtailments. 
However, they used a simple Linear Programming approach, which ignored discrete nature of operating modes. \citet{Ierapetritou2002} afterwards extended the approach to include discrete decisions using a mixed integer  model,  later improved by \citet{Karwan2007}. \citet{Zhang2015} investigated the economic benefit of energy storage in increasing the plant flexibility and the possibility of selling energy back to the market. In addition, they were among the first to propose a MILP model for participation of an energy-intensive process in the ancillary market. In order to guarantee the plan feasibility under all possible dispatching scenarios, they applied robust optimization technique. 
In spite of a high level of conservatism, their model proved to be profitable thanks to a large inventory connected to an on-site generation unit. This was demonstrated for a week-long plan with worst dispatching scenario of once a day considering spot and reserve price to be known for the whole planning horizon.\\
\citet{Vujanic2012} had previously developed a robust optimization for scheduling of a batch plant where the reserve is crudely provided through shifting the start time of jobs with an a priori known consumption profile. 
Later \citet{Zhang2016} improved the conservatism of a static robust optimization by applying the adjustable counterpart considering reactive actions based on dispatch realization. However, it could be challenging to maintain the tractability if complex decision rules are needed. In a different approach toward demand response in the air separation plant, \citet{Xenos2016} treated the problem of cost minimization and reserve provision sequentially in a deterministic fashion. It is worth mentioning that all these methods used a discrete-time presentation to formulate a MILP model. Instances of other processes for which DSM solution has been offered includes: cement \cite{Mitra2012}, desalination \cite{Ghobeity2010}, pulp \cite{Hadera2015} and chlorine plants \cite{Babu2008}.\\

This work concerns another extremely energy-intensive sector, iron and steel, the largest electricity intensive branch in Europe after chemical industry.\\ The first systematic study on the potentials of load shifting in steel plants was conducted by \citet{Ashok2006}. He formulated a MILP scheduling model based on discrete-time representation to reduce peak period demand and electricity cost in response to a time-of-use (TOU) tariff. \citet{Zhang2010} later incorporated in the model more characteristics of process and solved it using Lagrangian decomposition.\\ The continuous-time formulation of steel plant scheduling was first studied by \citet{Nolde2010}. They used the proposed model for electricity load tracking of a stainless-steel plant. Their model was reformulated by \citet{Hait2011} to achieve a better computational performance, and later extended by \citet{Hadera2013} to incorporate parallel machines at the production stage. Recently \citet{Hadera2015a} expanded it to include on-site generation and the possibility of selling electricity to the grid. In the same year, \citet{Castro2013} developed a discrete model based on resource-task network (RTN) concept for melt shop scheduling in steel-making process under energy constraint and investigated the response of their model to varying energy price. In addition, they studied the trade-off between the level of modelling accuracy and the computational expenses incurred by the required time granularity.\\
\txtred{The most limiting assumption common among all of the proposed continuous-time approaches is that the consumption profile remains constant over the time span of the task. In discrete RTN-based models, on the other hand, the piece-wise constant power profile is predefined.}
In addition to making it impossible to incorporate time-dependent process constraints, such an assumption results in a sub-optimal solution when the process involves parallel machines.\\
Other manufacturing processes in iron and steel sector have not received the same level attention as has the steel-making process. In one of the few examples, an energy scheduling scheme was introduced by \citet{Artigues2013} for a foundry with multiple induction furnaces. Considering a simplified energetic model of the process and human resource availability constraint, they proposed a bi-level constraint programming/MILP approach which iteratively solves sequencing and energy scheduling problem. However, they verified this model only from a feasibility point of view.\\

\subsection{Contribution and Organization of the paper}
In this paper, we propose a novel scheduling tool for the optimal operation of an energy-intensive plant considering both energetic and productive aspects, with emphasis on the interaction with electricity market. 
This tool is an enabler for the implementation of demand response strategies in industrial applications targeting both energy and reserve  markets. 
Each of these markets is characterized by different attribute of risk and access requirements, which are addressed in two separate scheduling modules: Energy-Aware Scheduling (EAS) and Reserve Scheduling (RS). Both are based on a MILP model that is formulated to encompass all the operating and safety conditions of the underlying process.
Furthermore, we propose various techniques to improve the computational tractability of the optimization problem at each stage considering the temporal limitation imposed by the market.\\
This paper illustrates the methodology on an exemplifying case of metal casting. To the best of our knowledge, no other previous work has provided such a comprehensive iDSM scheme for this application, taking into account all relevant aspects of the process dynamics and limitations. 
In particular, the novel reserve provision scheme developed in this work is proved to be an effective tool for exploiting the residual flexibility in the capacity market. \\
{The rest of this paper is organized as follows: Section \ref{sec:prob_desc} gives a general description of the problem, presenting all relevant aspects of energy and capacity markets and a more comprehensive view of the proposed Demand-Response framework. In Section \ref{sec:proc_mod} a discrete-time model of process is developed. The resulting MILP model takes into account all operational constraints of process within a metal casting foundry. In Section \ref{sec:DR_mod}, the demand response optimization is proposed for both the day-ahead spot and reserve markets. In Section \ref{sec:CS} a case-study illustrates market participation of a multi-line foundry in Nordic markets. Finally, the paper draws some conclusions and remarks on future directions.}%


\Nomenclature[Q, 02]{$ f \in \mathcal{F} $}{Melting furnaces}
\Nomenclature[Q, 01]{$ c \in \mathcal{C} $}{Casting furnaces/lines}
\Nomenclature[Q, 01]{$ d \in \mathcal{D} $}{Reserve bid slots}
\Nomenclature[Q, 05]{$ l \in \mathcal{L} $}{Power units}
\Nomenclature[Q, 06]{$ m \in \mathcal{M} $}{Melting jobs}		
\Nomenclature[Q, 03]{$ j \in \mathcal{J} $}{Melting furnace stages}	
\Nomenclature[Q, 04]{$ k \in \mathcal{K} $}{Global time-grid}
\Nomenclature[Q, 07]{$ q \in \mathcal{Q} $}{Baseline/Reserve grid}

\Nomenclature[S, 01]{$\mathcal{F}_{c}$}{Furnaces serving the casting line $c \in \mathcal{C}$}	
\Nomenclature[S, 02]{$\mathcal{F}_{l}$}{Furnaces powered by the power unit $l \in \mathcal{L}$}	
\Nomenclature[S, 03]{$\mathcal{J}_{\scriptscriptstyle m,f}$}{Stages in melt cycle $ m $ of furnace $ f $ }	
\Nomenclature[S, 04]{$\mathcal{J}^{\scriptscriptstyle E}$}{Energy-dependent stages }
\Nomenclature[S, 05]{$\mathcal{J}^{\scriptscriptstyle T}$}{Time-dependent stages }	
\Nomenclature[S, 06]{$\mathcal{K}_{\scriptscriptstyle q}$}{Grid points in baseline/reserve interval $q$ }	
\Nomenclature[S, 07]{$\mathcal{Q}_{\scriptscriptstyle d}$}{Reserve intervals $q$ in bid slot $d$  }

\Nomenclature[U, 06]{\footnotesize  tp}{Tapping}	
\Nomenclature[U, 01]{\footnotesize  0}{Initial}
\Nomenclature[U, 05]{\footnotesize  tl}{Tapping ladle}
\Nomenclature[U, 02]{\footnotesize  bid}{Bidding}
\Nomenclature[U, 02]{\footnotesize  bl}{Baseline}
\Nomenclature[U, 04]{\footnotesize  re}{Reserve market}
\Nomenclature[U, 03]{\footnotesize  da}{Day-ahead market}
\Nomenclature[U, 00]{\footnotesize  *}{Optimal}

\Nomenclature[D, 01]{$x_{\scriptscriptstyle f,m,j}^{\scriptscriptstyle k}$}{Stage Activation}
\Nomenclature[D, 02]{$y_{\scriptscriptstyle f,m,j}^{\scriptscriptstyle k}$} {Supplementary features \\(semi-continuity, pre-emption, power rate) }
\Nomenclature[D, 03]{$\nu_{\scriptscriptstyle q,c}$}{Positive imbalance direction}
\Nomenclature[D, 03]{$\zeta_{\scriptscriptstyle f,m,j}$}{Day-After Flexibility (DAF) for reserve}
\Nomenclature[D, 04]{$\beta_{\scriptscriptstyle d,c}, \, \mu_{\scriptscriptstyle d,c}$}{Multiple reserve activation }

\Nomenclature[C, 03]{$t_{\scriptscriptstyle f,m,j}$}{Stage starting time \nomunit{$[s]$}}
\Nomenclature[C, 01]{$p_{\scriptscriptstyle f,m,j}^{\scriptscriptstyle k}$} {Power \nomunit{$[kW]$}}
\Nomenclature[C, 04]{$R_{\scriptscriptstyle q,c}$}{Power Reserve \nomunit{$[kW]$}}
\Nomenclature[C, 02]{$p_{\scriptscriptstyle q,c}^{\scriptscriptstyle \rm{bl} +/ \rm{bl} -}$} {Positive/negative imbalance power \nomunit{$[kW]$}}

\Nomenclature[P, 20]{$ \alpha_{\scriptscriptstyle f,j} $}{Energy correction parameter \nomunit{$[-]$}}
\Nomenclature[P, 201]{$ \gamma_{\scriptscriptstyle c} $}{Poring furnace power coefficient \nomunit{$[kW/m^3]$}}
\Nomenclature[P, 21]{$ \lambda_{\scriptscriptstyle t}^{\scriptscriptstyle \rm{m}}$}{Price at time $t$ in market $\rm{m}$ \nomunit{$[$ \euro $/kWh]$}}				
\Nomenclature[P, 21]{$ \lambda_{\scriptscriptstyle t}^{\scriptscriptstyle \rm{re}}$}{Reward for reserve availability  \nomunit{$[$ \euro $/kW]$}}				
\Nomenclature[P, 0]{$\bar{k}_{\scriptscriptstyle c,n}$}{Breakpoint $n$ on casting rate for line $c$ \nomunit{$[s]$}}		
\Nomenclature[P, 03]{$ r^{\scriptscriptstyle \rm pwr}$}{Power rate limit \nomunit{$[kW/s]$}}		
\Nomenclature[P, 02]{$ r^{\scriptscriptstyle \rm of}$}{Overflow rate limit \nomunit{$[kW]$}}	
\Nomenclature[P, 04]{$ r^{\scriptscriptstyle \rm  sp}$}{Splash rate limit \nomunit{$[kW]$}}	
\Nomenclature[P, 01]{$ r^{\scriptscriptstyle \rm  cast}$}{Casting rate \nomunit{$[m^3/s]$}}			
\Nomenclature[P, 0910]{$v$}{Pouring furnace volume \nomunit{$[m^3]$}}	
\Nomenclature[P, 0911]{$\breve{v}$}{Tapped volume into pouring furnace \nomunit{$[m^3]$}}
\Nomenclature[P, 0912]{$\invbreve{v}$}{Cast volume from pouring furnace \nomunit{$[m^3]$}}	
\Nomenclature[P, 193]{$ L_{\scriptscriptstyle \max}$}{Maximum number of ladles \nomunit{$[-]$}}
\Nomenclature[P, 36]{$ \Pi$}{Maximum likelihood of reserve activation \nomunit{$[-]$}}				
\Nomenclature[P, 30]{$ \hat{\Delta}$}{Stage time duration (minimal)\nomunit{$[s]$}}	
\Nomenclature[P, 31]{$ \Delta_{f,c}^{\overleftrightarrow{}}$}{Ladle cycle \nomunit{$[s]$}}
\Nomenclature[P, 190]{$ C$}{Spot market objective function \nomunit{$[$ \euro $]$}} 
\Nomenclature[P, 191]{$ E\, , \, \hat{E}$}{Stage energy (actual, min.) \nomunit{$[kWh]$}} 
\Nomenclature[P, 22]{$ \bar{\tau}$}{Max. time without reheating \nomunit{$[s]$}} 
\Nomenclature[P, 192]{$\hat{E}^{\scriptscriptstyle \rm sp}$}{$\min {E}_{\scriptscriptstyle f,m,2}$ with risk of splash \nomunit{$[kWh]$}} 
\Nomenclature[P, 1921]{$N_{\scriptscriptstyle q}$}{Market discretization to grid step ratio \nomunit{$[-]$}} 
\Nomenclature[P, 32]{$\hat{\Delta}^{\scriptscriptstyle \rm of}$}{min. time with risk of overflow \nomunit{$[s]$}} 
\Nomenclature[P, 34]{$\Delta E_{\scriptscriptstyle c}$}{Energy shifted in DAF mode for line $c$ \nomunit{$[s]$}} 
\Nomenclature[P, 341]{$\Delta T_{\scriptscriptstyle \rm m}$}{Time discretization step for market $m$ \nomunit{$[s]$}} 
\Nomenclature[P, 35]{$\delta t$}{Time discretization step \nomunit{$[s]$}} 
\Nomenclature[P, 194]{$ P^{\scriptscriptstyle \max}$}{Maximum Power \nomunit{$[kW]$}} 
\Nomenclature[P, 195]{$ P^{\scriptscriptstyle \max}_{\scriptscriptstyle l}$}{Maximum Power of power pack $l$ \nomunit{$[kW]$}} 
\Nomenclature[P, 196]{$ P^{\scriptscriptstyle \min / \max}_{\scriptscriptstyle f,m,j}$}{Min/Max stage dependent power \nomunit{$[kW]$}}

\section{Problem description} \label{sec:prob_desc}
\txtred{This work addresses the iDSM of  electricity-intensive batch processes with a multi-stage multi-line production system, where each stage is defined by temporal or energetic targets and it is assumed to be modulable.
To provide Demand Response services, the process has to rely on its flexibility. A Mixed Integer Linear Programming model is therefore conceived to represent process load elasticity. 
In addition to the general flexibility of the stages, the model  has to be tailored to include all relevant case-specific process dynamics. 
In order to maximize the DR benefits via process scheduling, an industrial DSM framework is conceived to address the different aspects of demand response.
To this aim, the optimization problem requires to incorporate market-related aspects. Since both energy and reserve market are considered, two different, but interacting optimization problems are presented.\\  
In what follows we discuss the technical requirements for entering energy and reserve market as a demand-response provider and how different components of the proposed DSM framework interact with each other.} %

\subsection{Market}

In most of Europe, the electricity is traded in spot energy markets taking place one day before the actual delivery. Such timing makes them fit for the participation of industrial units with daily production plans. It is worth mentioning that though in this work we only consider DR as a buyer in the wholesale market, one may use the same model to sell the excess from a long-term energy contract, in which case it would be characterized as dispatchable demand-response. \\ 
In spot markets, the price is mostly set through marginal price auction with hourly bids. This price setting mechanism allows price-independent bids, which guarantee the purchase of critical loads at the system price without risking the rejection.  
Accepted bids in the operation day should be followed closely and any deviation from baseline is settled ex-post. There exist two balance settlement systems: i) one-price system, where load imbalance is settled at a single price regardless of whether it is aggravating or reducing the system overall imbalance and ii) two-price model, where different prices are assigned based on the direction of imbalance market. In order to encourage DR participation, some grid operators use the one-price model for load BRPs.\\
Unlike the spot market, reserve markets are relatively new and regulations vary a lot from one country to another. 
Although the majority of ancillary markets in Europe are considering bigger roles for demand flexibility as a reserve source, only in few Nordic countries daily auctions are held for bidding on reserve capacities, while other major markets only offer long-term contracts for their reserve products, which makes them not suitable for industrial sectors with non-identical daily demands.
Bid symmetry condition is another significant market barrier: it requires reserve capacity to be provided in both the up and down direction at the same time. Those TSOs that are tending to engage DR allow asymmetrical bids. We assume this condition to be true in this paper, as we consider only the case of positive reserve (up-regulation). However, even though there is little capacity to offer as the negative reserve in energy-intensive industries, the method can be easily extended to include such a provision if applicable and justifiable (economically and technically). \\
Depending on the target market, there are various technical requirements to be included in the DR model. This includes response time, sustain time, activation frequency, bid time granularity and bid size granularity. Primary markets have more restrictive requirements regarding the response time, while on the other hand, their maximum sustain time is much shorter than the secondary and tertiary market, which makes them more suitable for energy-constrained processes. In addition, due to a lower number of qualified providers, the rewards for primary reserves are significantly higher with respect to other markets. Consequently, even though we are not going to make any market-specific assumption during the development of DR model, primary reserve market will be the subject of the case study in this paper.\\ 
Price predictability is another deciding factor influencing the profitability of DR model. While many studies have suggested effective models to forecast the spot price for various energy markets \cite{weron2014electricity}, the behaviour of balancing markets is difficult to predict as they are ideally designed to deal with severe unforeseen events. Accordingly, a dispatchable DR model should be economically tuned against such uncertainty, particularly when pay-as-bid is used as the price-setting system.
In the process of developing DR strategies, we have tried to keep related models general in terms of market specifications such as clearing mechanisms, gate closure, and time intervals.  However, when it is necessary to make assumptions, we will adopt the rules of Danish market (as the most advanced example in Europe), which we will eventually use for the case study. \\
In addition, we are going to assume the DR model is merely a price taker and its bidding has no effect on the clearing price of the market.

\begin{figure*}[!htb]
	\centering
	\scriptsize
		\includegraphics[scale=.7]{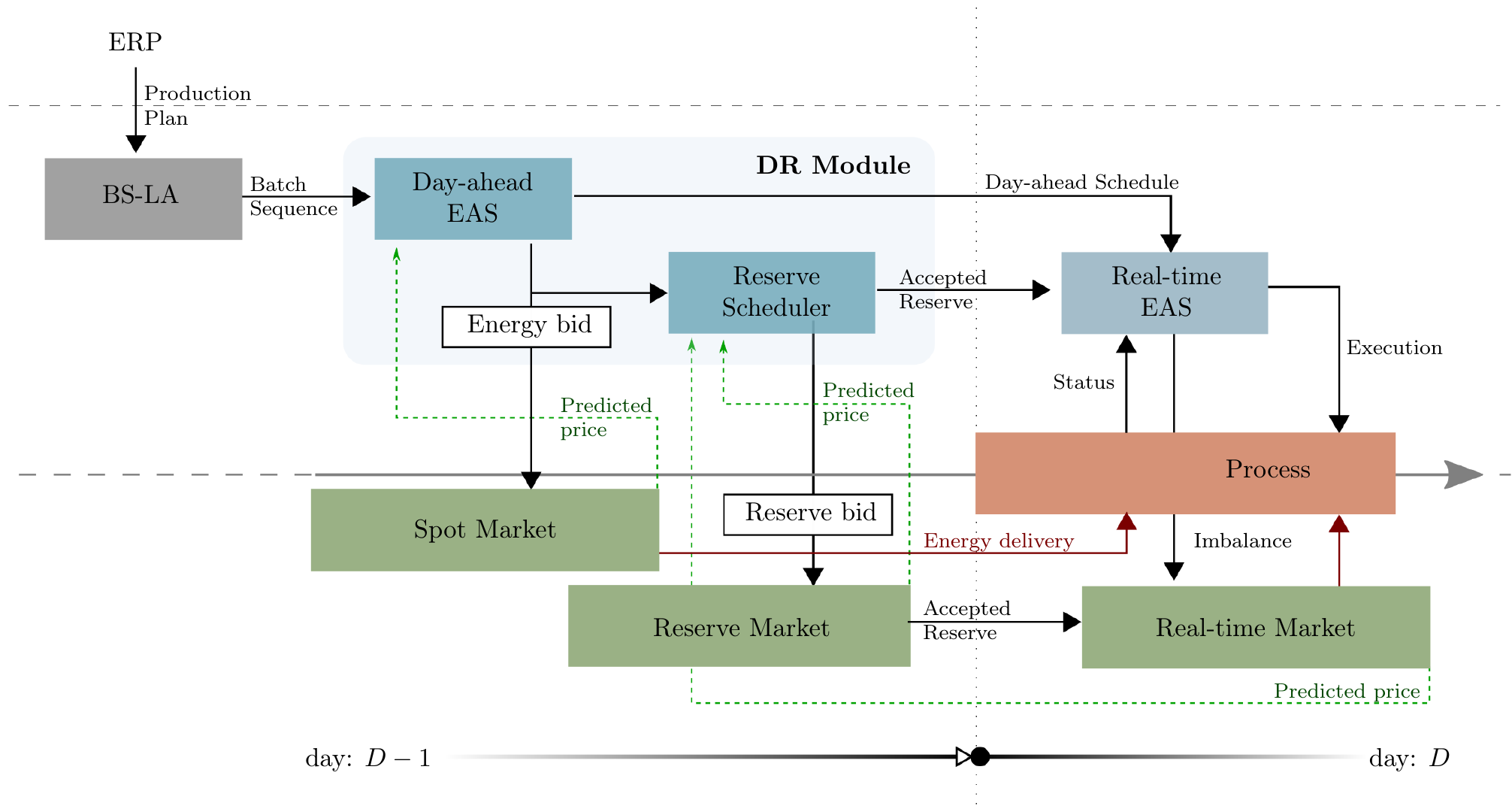}	
	\caption{Comprehensive framework to support industrial Demand Response. The scheme shows the interaction of the different components for iDSM (in the upper part) and the coupling with the corresponding electricity markets (in green), highlighting the time framework in which each one operates. The DR modules, which are the core of the paper, are framed in the grey box}
	\label{fig:frame}
\end{figure*} 

\subsection{Demand-response framework}
An ideal DR scheduling should be modelled such that the welfare is optimized simultaneously for both energy and reserve market. On the other hand, these two markets are very different in both their technical aspects and the rules governing them. In addition, the nature of these markets makes the risk associated with them vary significantly. Consequently, it would be very challenging to reach a monolithic model that could concurrently address different levels of priority on demand side for entering those markets.
The main goal of the iDSM tool in power-intensive sectors is to guarantee delivery of the energy required by the process to meet its delivery commitments. It is only when such a prerequisite is fulfilled that it may consider participating in a reserve market as a marginal opportunity. 
Consequently, similar to the work of \citet{Xenos2016}, we are going to formulate the participation in energy and reserve market separately. Firstly, we obtain a cost-minimizing power profile that  satisfies all operational and security constraints, then we calculate the maximum capacity that can be offered at each bidding interval assuming previously bid energy to be fully delivered.      

\Cref{fig:frame} illustrates the interactions of the DR module with the planning and scheduling and various electricity markets. 
\txtred{For a short-term scheduling, we distinguish two types of problems based on their nature and characteristic time scales: i) batch sequencing and line allocation (BS-LA), where the production plan given by Enterprise Resource Planning (ERP) is split into jobs, assigned to different lines, and  ii) DR module, subdivided in Energy-Aware Scheduler (EAS) and Reserve Scheduler (RS), for interacting with the spot and ancillary markets respectively.
Decisions taken in BS-LA can be generally excluded from the DR module. In case of interchangeable machines, sequencing and allocation have no significant effect on the optimality of the proposed scheme when it comes to energy cost minimization and reserve provision maximization. Due to different time-scales, such decoupling of the scheduling problems has little effect on the optimality of the resulting DR, while it drastically improves the computational tractability of the overall scheduling model.\\
The EAS layer receives these data and performs the intra-day scheduling, and it is  responsible for the short-term and detailed energetic aspects of the process. Using the optimal baseline, RS optimizes capacity provision considering potential imbalance costs.\\}
On the operation day, real-time EAS modifies the day-ahead scheduling based on the updates it receives from the process and considering the committed reserve. It minimizes the deviations from baseline and accordingly buys or sells the imbalance in the real-time market. The focus of this paper is on the modelling of DR module for participation in the day-ahead energy and reserve market. Development of the real-time module will be outside the scope of this paper and will be addressed by the authors in future work.

\section{Process Modelling}  \label{sec:proc_mod}
\txtred{The proposed DR optimization is based on a mixed integer model of the process flexibility. In this section we present a methodology to model both general and application-specific aspects of a metal casting process.}\\

\begin{figure}[!b]
	\centering
	\includegraphics[scale=0.8]{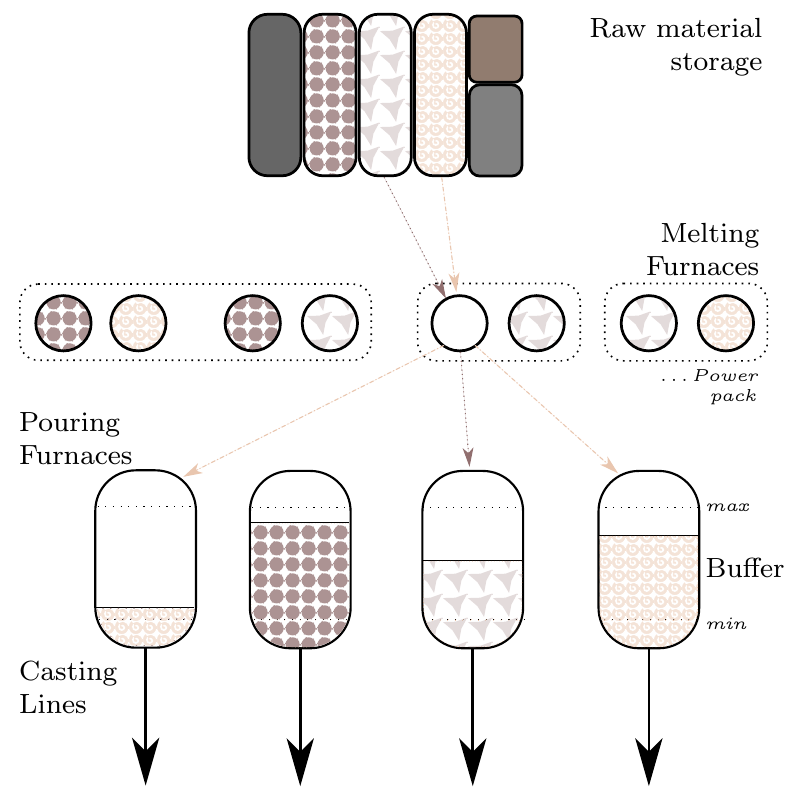}
	\caption{General work-flow in metal casting foundries. 
		Batch sequencing problem and line allocation alternatives are conceptually shown.
	}
	\label{fig:general_workflow}
\end{figure}
\begin{figure}[tb]
	\centering
	\includegraphics[scale=0.8]{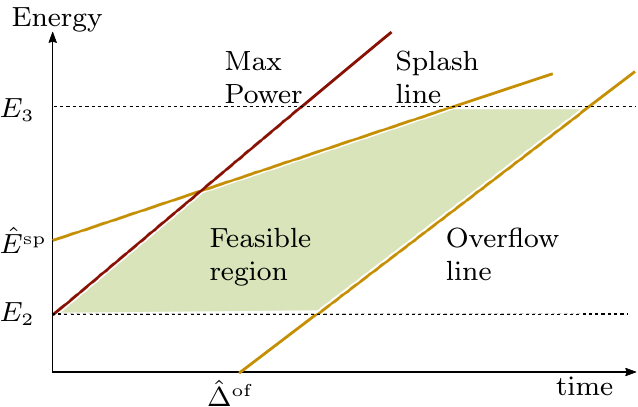}
		\caption{Feasible region of the Charge-Melting stage: heating and charging have to be interlinked to avoid the splash of overheated molten metal or the overflow of charged scrap}
	\label{fig:splash}	
\end{figure}

%
The different production stages within a multi-product foundry are illustrated in \Cref{fig:general_workflow}. Depending on capacity, a generic plant configuration includes a number of parallel casting lines, $c \in \mathcal{C}$, equipped with pouring furnaces that receive batches of molten metal from melting furnaces, $f \in \mathcal{F}$. Melting furnaces are fed with different types of scrap based on the alloy composition of the final product, defined in the production plan generated by the ERP system, considering the client orders. Based on the characteristics of the orders - i.e. material, number of products and their dimension - at each stage different decisions should be made regarding the unit assignment, lot sizing and sequencing of corresponding tasks.\\ 
Melting furnaces are grouped based on the common power line they may share. Accordingly,
we define $\mathcal{F}_l\subset \mathcal{F}$ as the subset of melting furnaces connected to a single power line, $l$. Pouring furnaces, on the other hand, are buffers  that receive batches of molten metal via ladles and fill sand moulds of different sizes at a fixed rate. The changeover of molten composition in each casting line poses a severe limitation on sequencing and lot sizing of the melt batches and their assignment to production lines.\\ Maximal productivity requires that the number of such changeovers is minimized and consequently the sequencing of  batches sharing the same material can be considered decision-making with a time-scale much longer than that of a day-long scheduling. \\
Accordingly, the Batch Sequencing and Line Allocation (BS-LA) defines the order of batches with the same material and assigns melting furnaces and casting lines such that the material availability and production rate set by ERP are respected. 
In general, this problem can be formulated as MILP with a simplified energy model similar to the work of \citet{Hadera2015}.   
The solution of BS-LA defines the number of melt cycles for each furnace, subset of furnaces $ \mathcal{F}_c \subset \mathcal{F}  $ assigned to a certain casting line $ c \in\mathcal{C} $ (such that  $\bigcap_{c\in \mathcal{C}}{\mathcal{F}_c}=\varnothing $). 
The DR module considers BS-LA input and the following model to optimize the energy provision.\\
\txtred{The proposed scheduling optimization is based on a discrete-time representation. The choice of this modelling approach, over the continuous one, can be justified mainly by its capability in integrating time-based events, such as complex process dynamics and  market coupling conditions. 
In addition, the solution optimality of a continuous model can be considered as an advantage only if the underlying process is deterministic. In a realistic case, the space left by the discretization-induced sub-optimality can be exploited to accommodate uncertainties of the process during the execution phase and avoid possible infeasibilities.
Accordingly, the planning horizon is discretized into a uniform grid, $\mathcal{K}$, and the power profile is piecewise constant and assumed controllable only on the grid nodes. The size of intervals, $\delta t$, should be chosen as the trade-off between the optimality and computational tractability and it should give an integer number of intervals within the day-ahead time slot, $\Delta T^{\scriptscriptstyle \text{da}}$, whose duration varies depending on the market rules.}
\begin{figure}[tb]
	\centering
	\includegraphics[scale=0.7]{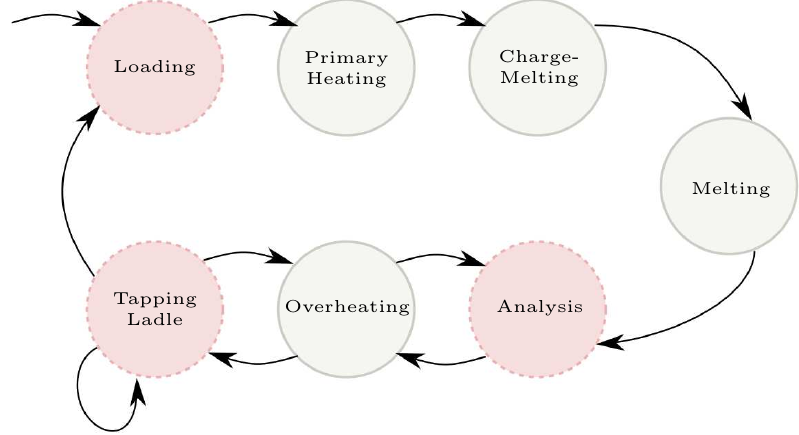}
	\caption{Finite State Machine of the melting furnace. Time-based (dotted contour) and Energy-based (solid) stages are shown. }
	\label{fig:FSM}
\end{figure}%
\txtred{\subsection{Melting operation}}\label{subsec:MeltOp}
Melting cycle is described by the finite-state machine in \cref{fig:FSM}. The stages are the following:
\begin{enumerate}[1)]
	\item \textit{Loading:} scrap metal is fed into the switched-off furnace. Time required for this stage sets a minimum interval between two melting jobs; \label{ph:A}
	\item \textit{Primary Heating:} the initial load is heated to form the molten bath, without introducing new scrap. Depending on the furnace technology, power ramp up rate is limited in this stage.\label{ph:B}
	\item \textit{Charge-Melting:} while heating continues, feeding resumes until the furnace capacity is full. The main challenge is synchronizing charging and melting so that there is neither a risk of splash (by overheating) nor furnace overflow (by overcharging); 	 
	\label{ph:C}
	\item \textit{Melting:} heating is continued until a target temperature is achieved; \label{ph:D}
	\item \textit{Analysis:} furnace is switched off for chemical analysis, composition correction and slag removal;  \label{ph:E}
	\item \textit{Overheating:} heating is resumed to raise the temperature to pouring point; \label{ph:F}
	\item \textit{Tapping:} molten metal is tapped into transportation ladles \txtred{multiple times}.  \label{ph:G}
\end{enumerate}
Phase transition is heterogeneously governed by exogenous and endogenous variables, such as temperature targets for stages \ref{ph:B}, \ref{ph:D} and \ref{ph:F}, weight for \ref{ph:A}, \ref{ph:C} and \ref{ph:G} and analysis approval for stage \ref{ph:E}. However, these transition variables can be effectively approximated by energy consumption and duration \txtred{of each stage}.
Accordingly, the process stages can be classified into two subsets, based on their completion requirements:
\begin{itemize}
	\item $ \mathcal{J}^{\scriptscriptstyle \text{E}} $ : Energy-based stages that are completed after consuming a certain amount of energy
	\item $ \mathcal{J}^{\scriptscriptstyle \text{T}} $ : Time-based stages that require no energy and whose completion is associated with a minimum duration
\end{itemize} 
The optimization variables of the EAS are the power values $ p_{\scriptscriptstyle f,m,j}^k $ defined at the grid nodes $k$ for each furnace $ f $, melt cycle $ m $ and stage $ j $ and the starting time of each stage, $t_{f,m,j}$. 
The proposed  modelling approach is based on the definition of boolean variable $X_{f,m,j}^k$ for each node on the time grid, \txtred{which is True for the nodes occurring after the start of stage $ j_{\scriptscriptstyle  f,m} $ and False for those preceding it}.\\
Considering that there is no waiting time between two consecutive stages, we can model timing and power consumption of each stage as the following logical statement:
\begin{equation}\label{eq:log_stat}
\txtred{\big[X_{\scriptscriptstyle f,m,j}^{\scriptscriptstyle k}\big] \Longleftrightarrow  \big[t_{\scriptscriptstyle f,m,j}\le  k\delta t\big]}  \qquad \forall f,m\in {\mathcal{M}_f},j,k
\end{equation}
\begin{multline}\label{eq:x_p}  
\big[\neg \left( X_{\scriptscriptstyle f,m,j}^{\scriptscriptstyle k} \wedge \neg X_{\scriptscriptstyle f,m,j+1}^{\scriptscriptstyle k} \right)\big]\implies \big[p_{\scriptscriptstyle f,m,j}^{\scriptscriptstyle k}=0\big]\\
\forall f,m\in {\mathcal{M}_f},j\in {\mathcal{J}^{\scriptscriptstyle \text{E}}} \setminus \left|\mathcal{J}\right|,k
\end{multline}  
While \txtred{the logical condition} \eqref{eq:log_stat} maps the $ t_{\scriptscriptstyle f,m,j} $ onto the node at which corresponding stage starts, the condition \eqref{eq:x_p} renders $p_{\scriptscriptstyle f,m,j}^{\scriptscriptstyle k}$, the power consumed by task $j_{\scriptscriptstyle f,m}$ zero at any instance out of the stage interval $\left[ t_{\scriptscriptstyle f,m,j},  t_{\scriptscriptstyle f,m,j+1} \right) $. \\
For each stage, we identify the minimum energy $ \hat{E}_{\scriptscriptstyle f,m,j} $  and corresponding time $ \hat{\Delta}_{\scriptscriptstyle f,m,j} $ leading to fastest completion of the stage.\\
However, in the modelling we need to consider the heat loss arising from any stage elongation imposed by energy-aware scheduling. Here we assume that a linear model can be fitted to give an over-estimation of heat loss in each stage as the function of its duration:
\begin{multline}   
E_{\scriptscriptstyle f,m,j}=\hat{E}_{\scriptscriptstyle f,m,j} \left(1+\sum\limits_{\scriptscriptstyle j''< j' \le j} \alpha_{\scriptscriptstyle f,m,j'}\frac{t_{\scriptscriptstyle f,m,j'+1}-t_{\scriptscriptstyle  f,m,j'}}{\hat{\Delta}_{\scriptscriptstyle f,m,j'}}\right) \qquad \\
\forall f,m\in \mathcal{M}_{\scriptscriptstyle f},j\in {\mathcal{J}^{\scriptscriptstyle \text{E}}}\label{eq:loss}
\end{multline}
where $\alpha_{ f,m,j}$ is the coefficient to be identified for each stage at minimum-time condition and $ j'' $ is the energy-based stage preceding $ j $ in $  j \in \mathcal{J}_E $.
Subsequently, the stage completion conditions for both energy-based and time-based stages give:
\begin{align}
&\sum\limits_{\scriptscriptstyle k\in \mathcal{K}}{p_{\scriptscriptstyle f,m,j}^{\scriptscriptstyle k}}\delta t = {E}_{\scriptscriptstyle f,m,j}&&\forall f,m\in \mathcal{M}_f,j\in \mathcal{J}^{\scriptscriptstyle \text{E}} \label{eq:Ed}\\ 
&t_{\scriptscriptstyle f,m,j+1}-t_{\scriptscriptstyle f,m,j}\ge \hat{\Delta}_{\scriptscriptstyle f,m,j}&&\forall f,m\in \mathcal{M}_{\scriptscriptstyle f},j\in {\mathcal{J}^{\scriptscriptstyle \text{T}}} \setminus \left|\mathcal{J}\right| \label{eq:Td} \\ 
&t_{\scriptscriptstyle f,m+1,1}-t_{\scriptscriptstyle f,m,\left|\mathcal{J}\right|}\ge \hat{\Delta}_{\scriptscriptstyle f,m\txtred{+1},\txtred{1}}&&\forall f,m\in \mathcal{M}_{\scriptscriptstyle f} \setminus \left|\mathcal{M}_{\scriptscriptstyle f}\right| \label{eq:Tdj}
\end{align}
In case the objective of the optimization involving constraint \eqref{eq:Ed} includes $ p_{\scriptscriptstyle f,m,j}^{\scriptscriptstyle k} $ for all $ k \in \mathcal{K} $ as minimization terms, the equality can be relaxed to inequality.\\
We assume that the end of tapping stage coincides with the start of loading stage for the next cycle, as any waiting time between two consecutive melting cycles can be included in the latter without limiting the feasible solution.\\
Using Big-M approach, \eqref{eq:log_stat} and \eqref{eq:x_p} can be formulated as mixed-integer linear inequalities. Subsequently, defining binary variable  $x_{\scriptscriptstyle f,m,j}^{\scriptscriptstyle k}$ with  correspondence to the Boolean $X_{\scriptscriptstyle f,m,j}^{\scriptscriptstyle k}$ the constraint \eqref{eq:log_stat} can be represented as:
\begin{multline}    
\txtred{\left(1 - x_{\scriptscriptstyle f,m,j}^{\scriptscriptstyle k}\ \right) \left(1+k \right) + K_{\scriptscriptstyle f,m,j}^{\scriptscriptstyle \min }} 
 \txtred{\le t_{\scriptscriptstyle f,m,j}/{ \delta t } } \\
 \txtred{\le x_{\scriptscriptstyle f,m,j}^{\scriptscriptstyle k}  + (1-x_{\scriptscriptstyle f,m,j}^{\scriptscriptstyle k}) K_{\scriptscriptstyle f,m,j}^{\scriptscriptstyle \max }} \qquad \forall f,m,j,k \label{xTd}
\end{multline}
\txtred{where $ K_{\scriptscriptstyle f,m,j}^{\scriptscriptstyle \min } $ and $ K_{\scriptscriptstyle f,m,j}^{\scriptscriptstyle \max } $  are the lower and upper bound on the time grid for each stage considering the minimum completion time of prior and later stages. Eventually these limits can be used to reduce the problem size. Further details can be found in \Cref{app:a1} } 
Similarly, constraint \eqref{eq:x_p} can be stated as:
\begin{multline}   
P_{\scriptscriptstyle f,m,j}^{\scriptscriptstyle \min }\left( x_{\scriptscriptstyle f,m,j}^{\scriptscriptstyle k}-x_{\scriptscriptstyle f,m,j+1}^{\scriptscriptstyle k} \right) 
\le p_{\scriptscriptstyle f,m,j}^{\scriptscriptstyle k} \\
\le P_{\scriptscriptstyle f,m,j}^{\scriptscriptstyle \max }\left( x_{\scriptscriptstyle f,m,j}^{\scriptscriptstyle k}-x_{\scriptscriptstyle f,m,j+1}^{\scriptscriptstyle k} \right) \\
\forall f, m \in {\mathcal{M}_f}, j\in {\mathcal{J}}\setminus \left| \mathcal{J} \right|, k \label{eq:xP}
\end{multline}
For each subset of furnaces connected to a single power line, we have the maximum power available limited by:
\begin{align}
\sum\limits_{\scriptscriptstyle f\in \mathcal{F}_l}{\sum\limits_{\scriptscriptstyle m\in \mathcal{M}_f}{\sum\limits_{\scriptscriptstyle j\in \txtred{\mathcal{J}}}{p_{\scriptscriptstyle f,m,j}^{\txtred{k}}}}}\le P_{\scriptscriptstyle l}^{\scriptscriptstyle \max }&&
\forall l,k  \label{Pmaxl}
\end{align}
Likewise, a limit on the overall power consumption is set by:
\begin{align}
\sum\limits_{\scriptscriptstyle f\in F}{\sum\limits_{\scriptscriptstyle m\in \mathcal{M}_f}{\sum\limits_{\scriptscriptstyle j\in \txtred{\mathcal{J}}}{p_{\scriptscriptstyle f,m,j}^{\scriptscriptstyle k}}}}\le {{P}^{\scriptscriptstyle \max }}&&\forall k  \label{Pmax}
\end{align}
By setting $ p_{\scriptscriptstyle f,m,j}^k\in \left\{ 0 \right\}\cup \left[ P^{\scriptscriptstyle \min }, \right.\left. P^{\scriptscriptstyle \max } \right] $, we can avoid a particular low-power  range without excluding the possibility for pre-emption (i.e. zero-power).\\
In order to make the power semi-continuous, we introduce binary variable  $y_{\scriptscriptstyle f,m,j}^{\scriptscriptstyle k}$ that enforce this feature using the constraint:
\begin{align}
P_{\scriptscriptstyle j}^{\scriptscriptstyle \min }y_{\scriptscriptstyle f,m,j}^{\scriptscriptstyle k}\le p_{\scriptscriptstyle f,m,j}^{\scriptscriptstyle k}\le P_{\scriptscriptstyle j}^{\scriptscriptstyle \max }y_{\scriptscriptstyle f,m,j}^{\scriptscriptstyle k}&&\forall f,m\in {\mathcal{M}_f},k
\end{align}
The same binary can be used to apply the power ramp limit, which is typical of stage \ref{ph:B}:
\begin{align}
p_{\scriptscriptstyle f,m,2}^{\scriptscriptstyle k}\le P^{\scriptscriptstyle \rm 0}+{{r}^{\scriptscriptstyle \rm pwr}}\sum\limits_{\scriptscriptstyle k'=1}^{\scriptscriptstyle k}{y_{\scriptscriptstyle f,m,2}^{\scriptscriptstyle k'}}\delta t&&\forall f,m \in {\mathcal{M}_f},k  \label{eq:p_rate_lim}
\end{align}
where $ P^{\scriptscriptstyle \rm 0} $ is the initial power 
and  ${r}^{\scriptscriptstyle \rm pwr}$ is the power rate limit. Inequality \eqref{eq:p_rate_lim} models a generalized power ramp limit, where the constraint increases only on time steps in which the power is on.\\
\txtred{\Cref{fig:splash} shows the feasible domain where the energy can evolve in stage \ref{ph:C} without the risk of splash or overflow.}
This region can be modelled using an MILP constraint:
\begin{multline}  
r^{\scriptscriptstyle \text{of}}\left( \sum\limits_{\scriptscriptstyle \txtred{k'=1}}^{\scriptscriptstyle \txtred{k}}{x_{\scriptscriptstyle f,m,3}^{\scriptscriptstyle \txtred{k'}}\delta t}-\hat{\Delta}^{\scriptscriptstyle \text{of}} \right)\le \sum\limits_{\scriptscriptstyle \txtred{k'=1}}^{\scriptscriptstyle \txtred{k}}{p_{\scriptscriptstyle f,m,3}^{{\scriptscriptstyle \txtred{k'}}}\delta t} \\
\le \hat{E}^{\scriptscriptstyle  \text{sp}}+r^{\scriptscriptstyle \text{sp}}\sum\limits_{\scriptscriptstyle \txtred{k'=1}}^{\scriptscriptstyle \txtred{k}}{x_{\scriptscriptstyle f,m,3}^{\scriptscriptstyle \txtred{k'}}}\delta t \qquad\forall f,m \in {\mathcal{M}_f},\txtred{k}
\end{multline}
$r^{\scriptscriptstyle \text{of}}$ and $r^{\scriptscriptstyle \text{sp}}$ are respectively the slopes of overflow and splash lines,  $\hat{E}^{\scriptscriptstyle \text{sp}}$ is the minimum energy at the end of stage \ref{ph:B} that could cause splash and  $\hat{\Delta}^{\scriptscriptstyle \text{of}}$ is the time it takes to fill up the furnace in absence of power, for a fixed rate of charge. All these quantities have to be identified from process data, for a nominal charge rate. 
\txtred{\subsection{Moulding operation }}
\txtred{The casting operation is concluded by the moulding process that includes heating and buffering of molten metal by pouring furnaces. 
The buffer is refilled by ladles, which deliver batches of molten iron from melting furnaces to each casting line. This transportation could be subjected to a limitation on ladle availability, which needs to be addressed.\\}

\subsubsection{\txtred{Casting buffer}}\label{subsec:PourF}
\txtred{The buffering elements of casting process are the pouring furnaces.} These furnaces act as the boundary between scheduling and planning layer in that they impose the pouring rate, with the holding volume constrained between an upper and lower limit. 
At each instant $ k \in \mathcal{K}$, we can impose these limits on the input-output balance of the molten iron from each pouring furnace as:
\begin{equation}
v_{\scriptscriptstyle c}^{\scriptscriptstyle \min } \le v_{\scriptscriptstyle c}^{\scriptscriptstyle 0} + \sum\nolimits_{f \in \mathcal{F}_c}{\breve{v}_{\scriptscriptstyle f}^{\scriptscriptstyle k-\Delta_{f,c}^{\overrightarrow{}}}} - \invbreve{v}_{\scriptscriptstyle c}^{\scriptscriptstyle k}   \le v_{\scriptscriptstyle c}^{\scriptscriptstyle \max} \label{eq:pouring_f_const}
\end{equation}
where $ v_{\scriptscriptstyle c}^{\scriptscriptstyle 0} $ is the volume of the buffer in the line $ c $, at the beginning of scheduling horizon, $ \breve{v}_{\scriptscriptstyle f}^{\scriptscriptstyle k} $ is the intake charge tapped from the melting furnace $ f $ and $ \invbreve{v}_{\scriptscriptstyle c}^{\scriptscriptstyle k} $  the throughput cast up to instant $ k $. $v_{\scriptscriptstyle c}^{\scriptscriptstyle \min } $ and $v_{\scriptscriptstyle c}^{\scriptscriptstyle \max } $ are the lower and upper limits of the buffer and $ \Delta_{\scriptscriptstyle f,c}^{\scriptscriptstyle \overrightarrow{}} $ is the overall time required for each ladle to deliver a batch of molten metal from melt cycle $ m_{\scriptscriptstyle f} $ to casting line $ c $.\\ 
\txtred{At a fixed production rate, defined as items per time unit, the rate at which the molten metal is poured from the buffer furnace can be altered, if the product type (e.g. mould dimension) changes.} Here we refer to the instances of product changeover and alteration of pouring rate as breakpoints, which together with the in-between rates define the continuous price-wise function for the cast volume in line $c$  up to the instance $k \in \left[ \bar{k}_{\scriptscriptstyle c,n},\bar{k}_{\scriptscriptstyle c,n+1} \right) $ as: 
\begin{equation}
\invbreve{v}_{\scriptscriptstyle c}^{\scriptscriptstyle k}= \sum\limits_{\scriptscriptstyle  n'=1}^{\scriptscriptstyle  N_c-1}{\bigg( \bar{k}_{\scriptscriptstyle c,n'} \left( {{r}_{\scriptscriptstyle c,n'+1}^{\scriptscriptstyle \rm cast}}-{{r}_{\scriptscriptstyle c,{n}'}^{\scriptscriptstyle \rm cast}} \right)-r_{\scriptscriptstyle c,n+1}^{\scriptscriptstyle \rm cast}k \bigg)}\delta t \label{eq:cast_metal}
\end{equation}
${{N}_{\scriptscriptstyle c}}$ is the number of intervals with different rates ${{r}_{\scriptscriptstyle c,n}^{\scriptscriptstyle \rm cast}}$ for each casting line $c$ and ${{k}_{\scriptscriptstyle c,n}}$ is the breakpoint on time grid indicating the beginning of interval $n$ on the line $c$.\\ 
\txtred{For a furnace with multiple tapping stage $ j \in \mathcal{J}^{\scriptscriptstyle \text{tp}} $, the molten volume delivered up to instant $ k $ can be calculated as:}
\begin{align}
\breve{v}_{\scriptscriptstyle f,m}^{\scriptscriptstyle k}=v^{\scriptscriptstyle \text{tl}}\sum\limits_{m\in \mathcal{M}_f}{\sum\limits_{j \in \mathcal{J}^{\scriptscriptstyle \text{tp}}}{{\txtred{x}}_{\scriptscriptstyle f,m,j}^{\scriptscriptstyle k} }}&&\forall f,m,k \label{eq:3tap}
\end{align} 
\txtred{where $ v^{\scriptscriptstyle \text{tl}} $ is the ladle effective volume.}\\  
Depending on furnace type, if the holding time takes longer than a minimum $ \bar{\tau}_{\scriptscriptstyle f,j}^{\scriptscriptstyle \text{tp}} $ 
a reheating should be performed, which can be enforced by:
\begin{multline} 
\sum\limits_{k\in \mathcal{K}}{p_{\scriptscriptstyle f,m,j}^{\scriptscriptstyle k}}\delta t\ge   \alpha_{\scriptscriptstyle f,m,j}\left( \left(t_{\scriptscriptstyle f,m,j+1}-t_{\scriptscriptstyle j}\right)/\bar{\tau}_{\scriptscriptstyle f,m,j}^{\scriptscriptstyle \text{tp}} \txtred{-1} \right)   \\
\forall f,m,j \in \mathcal{J}^{\scriptscriptstyle \text{tp}} \setminus \left|\mathcal{J}\right| 
\end{multline}
where $ \alpha_{\scriptscriptstyle f,m,j} $ is an experimentally identified coefficient similar to the one in equation (\ref{eq:loss}). \\
The power used by pouring furnaces to hold the temperature depends on their charged volume and here we assume that such a correlation can be over-approximated using a linear model, which allows us to estimate the total power consumed by holding furnaces at instant $ k $ as:
\begin{equation} \label{buffer_power}
\tilde{p}_{\scriptscriptstyle c}^{\scriptscriptstyle k}=\gamma_{\scriptscriptstyle c} \left(v_{\scriptscriptstyle c}^{\scriptscriptstyle 0}+\sum\limits_{f\in \mathcal{F}_c}{\sum\limits_{m\in \mathcal{M}_f}{\breve{v}_{\scriptscriptstyle f,m}^{\scriptscriptstyle k}-\invbreve{v}_{\scriptscriptstyle c}^{\scriptscriptstyle k}}}\right)
\end{equation}
where $ \gamma_{\scriptscriptstyle c} $ is a coefficient indicating power per volume required to maintain the molten in the desired temperature range.
\subsubsection{Transposition ladle}
\begin{figure}
	\centering
	\scriptsize
    \includegraphics[scale=1]{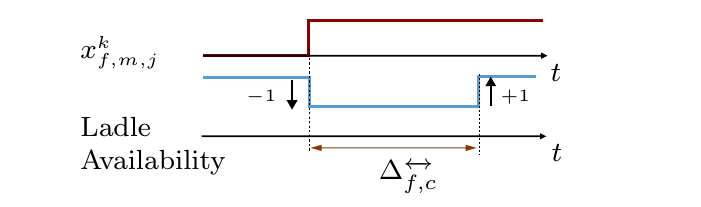}
	\caption{Binary representation of ladle availability: the variable ${x}_{\scriptscriptstyle f,m,j}^{\scriptscriptstyle k}$ is associated to tapping stages, i.e. $ \forall j \in \mathcal{J}^{\scriptscriptstyle \rm tp}$, while   $\Delta^{\leftrightarrow}_{f,c}$ is the time required by the ladle to reach the casting furnace, download the content and return.  }
		\label{fig:ladle}	
\end{figure}

In case the number of transportation ladles serving a casting line is limited, their scheduling must be modelled and linked to the timing of tapping. \\
As illustrated in Figure \ref{fig:ladle}, the change in ladle availability can be modelled via the linking binaries corresponding tapping stage. Once a tapping starts, the availability drops by one and its value is restored only after $ \Delta_{\scriptscriptstyle f,c}^{\overleftrightarrow{}} $, which is the overall time it takes for the ladle to deliver a batch from furnace $ f $ to casting line $ c $ and return to the same furnace array, including the loading and unloading time. Based on this concept, we can impose an upper bound of $ L^{\scriptscriptstyle \max} $  on the total number of ladles using constraint \eqref{eq:ladle}:
\begin{equation}
\sum\limits_{c\in C}\sum\limits_{f\in \mathcal{F}_c}\sum\limits_{m\in \mathcal{M}_f}\sum\limits_{j \in \mathcal{J}^{\scriptscriptstyle \text{tp}}}\left(\txtred{x}_{\scriptscriptstyle f,m,j}^{\scriptscriptstyle k}- \txtred{x}_{\scriptscriptstyle f,m,j}^{\scriptscriptstyle k-\Delta_{\scriptscriptstyle f,c}^{\overleftrightarrow{}}}\right) \le L^{\scriptscriptstyle \max}  \label{eq:ladle}  
\quad \forall k
\end{equation} %
\section{\txtred{Demand Response Optimization Problem } } \label{sec:DR_mod}
The flexibility of the energy-intensive process, mathematically modelled in the previous section, can be exploited within both implicit and explicit DR programmes. Each programme is formulated as a MILP with the general objective of maximizing the benefit of the industrial consumer. The appropriate objective function has to be particularized for the selected electricity market. 
All the requirements imposed by the energy or capacity market for such participation are defined as additional constraints. 

\subsection{Spot market}
The developed model allows the foundry to participate in the wholesale day-ahead energy market. 
We assume that the electricity price in day-ahead market has a time resolution of $\Delta T_{\rm{da}}$. However, the discretization at which the energy trade is balanced is typically a unit fraction of this time,  $\Delta T_{\rm{bl}}$. The latter is also the same sample time at which the baseline is set. Accordingly, the time interval of the scheduling grid, $\delta t$,  has to be chosen such that \txtred{ $ \Delta T_{\rm{bl}}$ is divisible by $\delta t $, i.e. }$\delta t | \Delta T_{\rm{bl}}$. \\
Considering the whole shop floor, the day-ahead optimization problem can be decomposed into sub-problems for each subset of furnaces $ \mathcal{F}_c $ exclusively assigned to casting line $ c $.  It must be assumed that a maximum power can be imposed separately for each subset.\\ 
The cost function for each casting line gives:
\begin{equation}
C_{\scriptscriptstyle c} = \sum\limits_{f\in \mathcal{F}_c} {\sum\limits_{m\in \mathcal{M}_{f}}}     {\sum\limits_{j\in \mathcal{J}}}   \sum\limits_{{k}\in {{\mathcal{K}}}}  \big[ {p_{\scriptscriptstyle f,m,j}^{\scriptscriptstyle {k}}} {\lambda^{\scriptscriptstyle \rm{da}}_{\scriptscriptstyle {k}}} \txtred{{\delta t}} \big] + \sum\limits_{k \in \mathcal{K}} \tilde{p}_{\scriptscriptstyle c}^{\scriptscriptstyle k}\lambda^{\scriptscriptstyle \rm{da}}_{\scriptscriptstyle k} \txtred{\delta t} \label{eq:spot_J}
\end{equation}
where $ \lambda^{\scriptscriptstyle \rm{da}}_{\scriptscriptstyle k} $ is the day-ahead market price.\\
Subsequently, for each line, solving the optimization problem:
\begin{align}
\txtred{\min_{p, x, y}} & \txtred{\quad C_{\scriptscriptstyle c}} \label{eq:spot_op} \\ 
\txtred{\text{s.t.} }&\txtred{\quad \text{eq.} \eqref{xTd}-\eqref{eq:cast_metal}}  \nonumber 
\end{align}
we obtain the optimal power profiles $ p_{\scriptscriptstyle f,m,j}^{* \scriptscriptstyle k  } $ that results in the minimum cost $ C_{\scriptscriptstyle c}^* $. \\
Consequently, the optimal baseline for each interval $ q \in \mathcal{Q}$ is calculated as:
\begin{align}
P_{\scriptscriptstyle q,c}^{\scriptscriptstyle\text{b}}= \frac{1}{N_q} \sum\limits_{k\in \mathcal{K}_q}\sum\limits_{f\in \mathcal{F}_c} \sum\limits_{m\in \mathcal{M}_f}\sum\limits_{j\in \mathcal{J}}{p_{\scriptscriptstyle f,m,j}^{ * \scriptscriptstyle  \txtred{k}}}  &&\forall q,c \label{Pbaseline}
\end{align}
where  \txtred{$\mathcal{K}_q = \{k \in  \mathcal{K} |  (q-1)N_{\scriptscriptstyle q}<k\leq qN_{\scriptscriptstyle q}\}$ 	is the subset of  $ \mathcal{K}$ corresponding to the baseline interval} and $N_{\scriptscriptstyle q} ={\Delta T_{\scriptscriptstyle \rm{bl}}}/{\delta t} $.\\
Eventually the total cost of energy consumed by the whole plant is given by
\begin{equation*}
C = \sum\nolimits_{c \in \mathcal{C}} C_{\scriptscriptstyle c}
\end{equation*} %
\txtred{It is worth mentioning that the additional fee regarding the maximum contracted power demand (MCD) have not been included in the cost function, as it is assumed to be negotiated for a much longer time period and thus constant during the scheduling horizon. In addition, we assume that the total installed capacity of all power units is equal to the MCD.}

\subsection{Reserve provision}
Participation in the ancillary market is essentially a matter of how much capacity can be offered without risking the infeasibility of a re-scheduling problem - if activated - and at what price such capacity should be offered to cover all the incurred expenses.
Although there is no unified regulation across Europe for capacity market, it is safe to assume the bidding interval to be an integer multiple of baseline time discretization and  the maximum duration of a single dispatch not to be less than this latter. 
Without loss of generality, we assume the reserve grid and the baseline to share one resolution.\\
Depending on response time requirements and the maximum duration, which varies considerably between different reserve markets, the interruption period may span across several time intervals of the grid.\\
Consequently, \txtred{the optimal capacity bid can be formulated as $ |\mathcal{Q}|$ separated multi-objective optimizations}, whose main goal is to maximize the uniform reduction of baseline $ P_{\scriptscriptstyle q}^{\scriptscriptstyle \rm{bl}} $ (as the maximum total power) for all $ q \in \left\{q,\dots,q+M\right\} $, with $ M $ being the maximum number of intervals that may be affected by reserve activation. The secondary objective of this optimization problem is minimizing the cost of the internal imbalance sustained as the result of reserve activation, since 
a definite amount of energy is necessary to fulfil the planned production.\\
While the offered capacity is weighted by the  reward for availability, which is guaranteed if bid accepted, the incurred cost of imbalance must be included considering the probability of the reserve being called upon in a specific time slot. \\   
Using the baseline computed for each casting line, $c$, in the spot-market problem, we set up a multi-objective optimization to calculate the optimal amount and price of reserve capacity bidding for all intervals $ q \in \mathcal{Q}_{\scriptscriptstyle d} $ within in the time block $ d $:
{\allowdisplaybreaks
\begin{IEEEeqnarray}{lcll} \label{explicitq} 
		\IEEEyesnumber  \IEEEyessubnumber*
	\max
	&\quad& \hspace{-5pt} \lambda_{\scriptscriptstyle q}^{\scriptscriptstyle \text{re}}R_{\scriptscriptstyle q,c} - \Pi \left( \sum\limits_{\scriptscriptstyle q'>q''} \!\! {\left(\lambda_{\scriptscriptstyle q'}^{\scriptscriptstyle +} p^{\scriptscriptstyle \text{bl}+}_{\scriptscriptstyle q',c} + \lambda_{\scriptscriptstyle q'}^{\scriptscriptstyle -} p^{\scriptscriptstyle \text{bl}-}_{\scriptscriptstyle q',c}\right)}\Delta T_{\scriptscriptstyle \rm{bl}}- \bar{\lambda} \Delta E_{\scriptscriptstyle c}\right)\nonumber\\
	\text{s.t.}
	&&\hspace{-15pt}  \sum\limits_{k\in \mathcal{K}^M_q}\sum\limits_{f\in \mathcal{F}_c} \sum\limits_{m\in \mathcal{M}_f}\sum\limits_{j\in \mathcal{J}}\frac{p_{\scriptscriptstyle f,m.j}^{\scriptscriptstyle \txtred{k}}}{N_q} \!\! \le P^{\scriptscriptstyle \text{bl}}_{\scriptscriptstyle q,c} -R_{\scriptscriptstyle q,c}    \label{exp1} 	\\	
	&&\hspace{-15pt}  \sum\limits_{k\in \mathcal{K}_{q'}}\sum\limits_{f\in \mathcal{F}_c} \sum\limits_{m\in \mathcal{M}_f}\sum\limits_{j\in \mathcal{J}}\frac{p_{\scriptscriptstyle f,m.j}^{\scriptscriptstyle \txtred{k}}}{N_q} = P^{\scriptscriptstyle \text{bl}}_{\scriptscriptstyle q',c}+ p^{\scriptscriptstyle \text{bl}-}_{\scriptscriptstyle q',c}+ p^{\scriptscriptstyle \text{bl}+}_{\scriptscriptstyle q',c}  \label{exp2} \\
	&&\hspace{-15pt}  0 \le  p^{\scriptscriptstyle \text{bl}+}_{\scriptscriptstyle q',c} \le \nu_{\scriptscriptstyle q',c}  P^{\scriptscriptstyle \text{bl}}_{\scriptscriptstyle q',c} \label{exp3}  \\
	&&\hspace{-15pt}  (\nu_{\scriptscriptstyle q',c}-1)   P^{\scriptscriptstyle \text{bl}}_{\scriptscriptstyle q',c}\le  p^{\scriptscriptstyle \text{bl}-}_{\scriptscriptstyle q',c} \le 0  \label{exp4} \\     
	&&\hspace{-15pt}  0 \le R_{\scriptscriptstyle q,c} \le P^{\scriptscriptstyle \text{bl}}_{\scriptscriptstyle q,c} \label{exp5}   \\
	\IEEEeqnarraymulticol{3}{r}{\forall q'>q''}\nonumber   \\ 
	&&\hspace{-15pt}  \text{and constraints } \eqref{eq:Ed}  \text{ to } \eqref{buffer_power} \nonumber  
\end{IEEEeqnarray}}%
where  \txtred{$\mathcal{K}^M_q = \{k \in  \mathcal{K} |  (q-1)N_{\scriptscriptstyle q}<k\leq q''N_{\scriptscriptstyle q}\}$ and $ q'' = \min(q+M, \max(\mathcal{Q}_d))$.}\\
In the objective function, $ R_{\scriptscriptstyle q,c} $ is the reserve power weighted by $ \lambda_{\scriptscriptstyle q,c}^{\scriptscriptstyle \text{re}} $, i.e. the estimated standby payment for upward regulation, while $ \Pi $ is the maximum likelihood of activation in the specific day which is used to scale all the incurred costs. These costs include, $ p^{\scriptscriptstyle \text{bl}+}_{\scriptscriptstyle q',c}/p^{\scriptscriptstyle \text{bl}-}_{\scriptscriptstyle q',c} $ , positive/negative power weighted by up/down regulation price $ \lambda_{\scriptscriptstyle q}^{\scriptscriptstyle +}/\lambda_{\scriptscriptstyle q}^{\scriptscriptstyle -} $. 
Constraint \eqref{exp1} imposes the capacity reduction and constraints \eqref{exp2}-\eqref{exp4} are in place $\forall q'>q''$ to calculate positive and negative power deviations. A binary variable $ \nu $  differentiates positive values from negative ones. Finally, constraint \eqref{exp5} sets bounds of the offered capacity.\\
In the case imbalances are to be settled using a one-price model (i.e. $ \lambda_{\scriptscriptstyle q}^{\scriptscriptstyle +} = \lambda_{\scriptscriptstyle q}^{\scriptscriptstyle -} $), the optimization can be simplified by using  the absolute value of deviation $ P_{\scriptscriptstyle bl}^{\scriptscriptstyle b\pm} $, which renders constraints \eqref{exp3} and \eqref{exp4} useless.   
\txtred{In case the dispatched energy is enumerated by TSO, then an amount of $ \Pi \lambda_{\scriptscriptstyle q}^{\scriptscriptstyle +} $ can be added to the standby price in the objective function.} However, we consider this increment to be negligible with respect to imbalance cost.\\
By solving the optimization \eqref{explicitq} for all intervals $ q \in \mathcal{Q}_{\scriptscriptstyle {d}} $ the optimal capacity that can be offered by each casting line $ c $  for bidding block $ d $ gives:
\begin{equation}
R^{\rm{bid}}_{\scriptscriptstyle d,c} = \min \left\{R^{\scriptscriptstyle *}_{\scriptscriptstyle q,c}\right\}_{q \in \mathcal{Q}_d} \label{r_bid}
\end{equation}
with the minimum bidding price of:
\begin{equation}
\lambda^{\rm{bid}}_{\scriptscriptstyle d,c} = \max \left\{\lambda^{\scriptscriptstyle *}_{\scriptscriptstyle q,c}\right\}_{q \in \mathcal{Q}_d}   \label{p_bid}
\end{equation}
where $ R^{\scriptscriptstyle *}_{\scriptscriptstyle q,c} $ is the optimal solution of the optimization (\ref{explicitq}) solved at instance $ q $ with the imbalance cost of  $ \lambda^{\scriptscriptstyle *}_{\scriptscriptstyle q,c} $ per capacity unit.\\
\txtred{It should be emphasised that each optimization generates, not only bidding quantities \eqref{r_bid}-\eqref{p_bid}, but also a contingency plan for each possible activation, as the updated optimal schedule.}\\

Depending on parallelization capability of the computation unit, optimization \eqref{explicitq} can be solved concurrently for more than one interval at a time. Although computational time of each optimization can be considerably reduced by warm-starting with the corresponding portion of the solution from  EAS model (i.e. the optimization starts from a feasible solution and search to improve it till its terminated due time constraints), it may be still a burden to solve all possible reserve optimization problems for all the bidding intervals of the day. \\
On one hand, due to low activation probability, the primary goal of bidding should be selling as much  capacity as possible, on the other hand, a large bid is more probable to be rejected due to congestion problems. The activation frequency is proportional to the share of accepted capacity with respect to the whole reserve market.\\
Whether a bid is accepted or not depends on  many factors including market capacity and acceptance mechanism. However, a quantitative evaluation of participation profitability is not possible without  modelling the entire market, which is outside the scope of this work.\\%
So far, we have considered the reserve optimization problem to have a one-way interaction with the spot market optimization through the baseline profile. However, it is possible to significantly increase revenue in the former with no or slight deterioration of optimal cost in the latter with a baseline closer to the maximum power limit. It is particularly beneficial for reserve biddings with high time granularity where a reserve capacity has to be sustained for a very long period, with respect to the baseline interval, in order to be accepted.\\ 
To address this limitation, we can add to the spot optimization problem \eqref{eq:spot_op} the following constraint:
\begin{align}
\underbar{P}_{\scriptscriptstyle d}^{\scriptscriptstyle \text{re}} \le P_{\scriptscriptstyle q,c}^{\scriptscriptstyle \text{bl} } && \forall d, q \in \mathcal{Q}_d
\end{align}    
introducing new variable $ \underbar{P}_{\scriptscriptstyle d}^{\scriptscriptstyle \text{r}} $ which is to be included in the objective function \eqref{eq:spot_J}, with a negative sign. Though as we have mentioned before, the overall priority should be spot market optimization, the weight associated variable $ \underbar{P}_{\scriptscriptstyle d}^{\scriptscriptstyle \text{r}} $  can be selected as a trade-off on how much of the flexibility is to be liquidized in the spot market and how much in the ancillary market.
\txtred{\subsubsection{Enhanced reserve provision : day-after flexibility}}
In addition to the flexibility from day ($ D $), a low probability of reserve activation may economically justify offering a portion of next day ($ D+1 $) flexibility to increase the reserve offered in the market. 
Here we model this day-after flexibility (DAF) and quantify the risk associated with any subsequent limitation imposed to the next day.
DAF implies a partial or complete shift for some of the last melt cycles to the next day. \\ 
To this end, we introduce the binary decision variable $ \zeta_{\scriptscriptstyle f,m,j} $, which if true implies that phase $ j_{\scriptscriptstyle f,m}  $ is scheduled to start on the same day and if false means that it has been moved to the next day.\\
Such a decision can be enforced by the following constraints for all stages set by planner for the current day:
\begin{multline}
\zeta_{\scriptscriptstyle f,m,j+1} \le  \left({\sum\nolimits_{\scriptscriptstyle k\in {\mathcal{\txtred{K}}}}{p_{\scriptscriptstyle f,m,j}^{\scriptscriptstyle k}}\delta t}\right) \left/ {\hat{E}_{\scriptscriptstyle f,m,j} } \le \zeta_{\scriptscriptstyle f,m,j} \right.\\ \forall f,m \in \tilde{\mathcal{M}}_{\scriptscriptstyle f},j\in {\mathcal{J}^{\scriptscriptstyle \text{E}}} \setminus \left|\mathcal{J}\right| \label{E2next}
\end{multline}
\begin{multline}
\zeta_{\scriptscriptstyle f,m,j+1} \le \left({t_{\scriptscriptstyle f,m,j+1}-t_{\scriptscriptstyle f,m,j}}\right)/{{\Delta}_{\scriptscriptstyle f,m,j}} \le \zeta_{\scriptscriptstyle f,m,j} \\ \forall f,m \in \tilde{\mathcal{M}}_{\scriptscriptstyle f},j \in {\mathcal{J}^{\scriptscriptstyle \text{T}}} \setminus \left|\mathcal{J}\right| \label{T2next}
\end{multline}
Where $ \tilde{\mathcal{M}}_f $ is the subset of relaxed cycles for each furnace. Decision on the extension of these subsets should be taken considering the computational cost due to model size increase and whether it can justify the possible benefits.\\
Similar to the spot market problem, each casting line is considered to be independent from the others at the scheduling level, which allows us to calculate the offered capacity for each one separately.\\
In order to avoid scheduling infeasibility due to buffer initial volume in the day after, we add a margin $v_{\scriptscriptstyle c}^{\scriptscriptstyle \text{s}}$ to the lower bound of constraint \txtred{\eqref{eq:pouring_f_const}} as: 
\begin{equation}
v_{\scriptscriptstyle c}^{\scriptscriptstyle \min } + v_{\scriptscriptstyle c}^{\scriptscriptstyle \text{s}} \le v_{\scriptscriptstyle c}^{\scriptscriptstyle 0} + \sum\nolimits_{f_c}{\breve{v}_{\scriptscriptstyle f}^{\scriptscriptstyle k-\Delta_{f,c}^{\overrightarrow{}}}} - \invbreve{v}_{\scriptscriptstyle c}^{\scriptscriptstyle k} \le v_{\scriptscriptstyle c}^{\scriptscriptstyle \max} \label{eq:buffer}
\end{equation}
For each furnace, this new lower bound should take into account the shifted time of relaxed melt cycles and add a proper margin to compensate the delay in delivery to the assigned casting line.
Thus for each line $ c $, assuming the casting rate to remain at the value of $ r_{\scriptscriptstyle c,N_c-1}^{\scriptscriptstyle \text{cast}} $,  the safety margin $ v_{\scriptscriptstyle c}^{\scriptscriptstyle \text{s}} $ can be calculated as:
\begin{equation}
v_{\scriptscriptstyle c}^{\scriptscriptstyle \text{s}} = r_{\scriptscriptstyle c,N_c-1}^{\scriptscriptstyle \text{cast}}  \sum\limits_{f\in \mathcal{F}_c}\sum\limits_{m\in \mathcal{M}_f}\sum\limits_{j \in \mathcal{J}} \Delta_{\scriptscriptstyle f,m,j}^{*} \left(1-\zeta_{\scriptscriptstyle f,m,j}\right) \label{eq:v_cs}
\end{equation} 
where $ \Delta_{\scriptscriptstyle f,m,j}^{*} $ is the duration of stage $ j_{\scriptscriptstyle f} $ of furnace $ f $ as scheduled by the day-ahead optimization. 
In order to quantify the associated risk, we calculate the amount of energy that has been rescheduled to be bought in the next day as:
\begin{equation}
\Delta E_{\scriptscriptstyle c} = \sum\limits_{f\in \mathcal{F}_c}\sum\limits_{m\in \mathcal{M}_f}\sum\limits_{j \in \mathcal{J}} \hat{E}_{\scriptscriptstyle f,m,j}\left(1-\zeta_{\scriptscriptstyle f,m,j}\right) 
\label{eq:DeltaE_c}
\end{equation} 
for all furnaces assigned to buffer $ c $. \\Depending on the estimation of process flexibility on the next day, this value can be limited by an upper bound to reduce the risk of feasibility, in addition to the safety margin set by \eqref{eq:v_cs}.\\ 
In order to quantify the incurred cost we assume that in case of reserve activation, the shifted energy is purchased at the maximum spot price of day $ D+1 $ and thus we include $ \Delta E_{\scriptscriptstyle c} $ in the objective function of explicit model, weighting it accordingly.     
The decision on whether to include the relaxation model for a bidding interval should be made considering the market volume and the portion such relaxation may contribute the total offered capacity. For instance, during early hours of scheduling the original flexibility obtained from the non-relaxed problem could provide enough capacity.
\txtred{\subsubsection{Multiple daily reserve activation : internal aggregation}}
A crucial factor in the reserve provision is the frequency by which the reserve is called upon. This not only affects the optimal capacity for each line when weighted in the multi-objective optimization, but also imposes a severe limitation on the exploitable capacity without resorting to aggregation.
If activation number is limited to once a day then the total capacity will be simply equal to $ \sum\nolimits_{\scriptscriptstyle c \in \mathcal{C}} R^{\scriptscriptstyle \rm{bid}}_{\scriptscriptstyle c,d} $ \txtred{minimally} priced at $ \sum\nolimits_{\scriptscriptstyle c \in \mathcal{C}} \lambda_{\scriptscriptstyle  c,d} R^{\scriptscriptstyle \rm{bid}}_{\scriptscriptstyle c,d} /\sum\nolimits_{\scriptscriptstyle c \in \mathcal{C}} R^{\scriptscriptstyle \rm{bid}}_{\scriptscriptstyle c,d} $.\\ \txtred{The same optimal solution cannot be obtained if the TSO request multiple activations per day, in which case, either internal or external should be considered for entering the ancillary market.}\\
In the internal aggregation we distribute the risk of dispatching among all the casting lines in the factory. Consequently, we can formulate the aggregation problem as the combinatorial optimization:
{\allowdisplaybreaks
\begin{IEEEeqnarray}{lcll} \label{multiact}
	\IEEEyesnumber  \IEEEyessubnumber*
	\max_{\beta,\tilde{R}} &\quad&  \sum\limits_{\scriptscriptstyle d\in \mathcal{D}} \left( \tilde{R}_{\scriptscriptstyle d} +\sum\limits_{\scriptscriptstyle c\in \mathcal{C}} \beta_{\scriptscriptstyle d,c} R^{\scriptscriptstyle \rm{bid}}_{\scriptscriptstyle d,c} \lambda^{\scriptscriptstyle\rm{bid}}_{\scriptscriptstyle d,c}\right) & \nonumber\\
	\text{s.t.}
	&& \sum\limits_{\scriptscriptstyle d\in \mathcal{D}} \beta_{\scriptscriptstyle d,c} \le 1 & \forall c  \\ \label{multiact1}
	&& \sum\limits_{\scriptscriptstyle c\in \mathcal{C}} \beta_{\scriptscriptstyle d,c} \le \left|\mathcal{C}\right| - N_{\scriptscriptstyle  a} + 1 &  \forall d  \\ \label{multiact2}
	&& \tilde{R}_{\scriptscriptstyle d} \le \left(1-\beta_{\scriptscriptstyle d,c}\right) R^{\rm \scriptscriptstyle bid, \max}_{\scriptscriptstyle d} & \forall d,c \\  \label{multiact3}
	&& \tilde{R}_{\scriptscriptstyle d} \le R^{\scriptscriptstyle \rm{bid}}_{\scriptscriptstyle d,c} + \left(1-\mu_{\scriptscriptstyle d,c}\right) R^{\scriptscriptstyle \rm{bid}, \max}_{\scriptscriptstyle d} \quad & \forall d,c\\  \label{multiact4}
	&& \sum\limits_{\scriptscriptstyle c\in \mathcal{C}} \mu_{\scriptscriptstyle d,c} \ge N_{\scriptscriptstyle a}-1 & \forall d	\label{multiact5}
\end{IEEEeqnarray}}
where $ N_{\scriptscriptstyle a} $ is the maximum number of dispatching per day and $ \beta_{\scriptscriptstyle d,c} $ is a binary variable which is equal to one if the capacity of line $ c $  is included in the total sum offered as one of the $ N_{\scriptscriptstyle a} $ largest reserve blocks and zero otherwise. Constraints \eqref{multiact1} and \eqref{multiact2} ensure that each block is considered just once in that capacity, while constraints \eqref{multiact3} to \eqref{multiact5} using the binary variable $ \mu $ are implemented to allocate $ N_{\scriptscriptstyle a}$-th minimum block if no other block has been picked from time interval $ d $. By including $ \tilde{R} $ we assume that reserve block can be partially dispatched if enforced and that price per units remains constant for all such partial realizations.\\
In case the contract with TSO/BRP includes a minimum period $ \Delta T^{\scriptscriptstyle \rm{bid}} > \delta q $ between two consecutive activation, the aggregated reserve can be improved by temporal decomposition of optimization \eqref{multiact} such that each problem is solved over $ \Delta T^{\scriptscriptstyle \rm{bid}}  $. \txtred{In case the number of possible activations requested bt TSO is more than $ \left|\mathcal{C}\right| $, external aggregation will be inevitable for participation in the reserve market.}\\
\section{Case Study and Results}  \label{sec:CS}
In this section we apply the proposed demand response optimization approach to a real-world manufacturing plant with melting, casting and machining departments.   
The foundry is a 300.000 ton/year plant with eight medium frequency coreless induction furnaces for the melting stage. 
Every couple of induction furnaces is connected to a single power unit (PU). 
This battery of four couples of induction furnaces delivers a charge of molten metal to the channel-type induction holding furnaces for the casting process, through ladles moved by forklifts. 
Due to confidentiality reasons we are not going to disclose the details of demand and exact consumption for this plant.\\ 
\begin{table}[!b]
	\caption{Time discretization steps}
	\label{table_example}
	\begin{center}
		\begin{tabular}{|l||l||l|l|}
			\hline
			Type & Step & Index & Value\\
			\hline
			Scheduling Optim. & $\delta t$ & $k$ & $5$ min\\
			\hline
			Baseline & $\Delta T_{\rm{bl}}$ & $q$ & $15$ min\\
			\hline
			Day-ahead Market & $\Delta T_{\rm{da}}$ & $-$ & $1$ h\\
			\hline
			Reserve Market & $\Delta T_{\rm{rm}}$ & $d$ & $4$ h\\
			\hline
		\end{tabular}
	\end{center}
\end{table}%
In the absence of a specific market model, whose development is outside the scope of this paper, in order to verify benefits of the proposed DR programme we have used the historical data from the western Denmark (DK1) bidding zone. The choice of this market is motivated by the leading-edge status of Nordic market, which allows a comprehensive analysis of participation in both the energy and capacity market.%

\subsection{Spot Market}
\begin{figure}[!t]
	\centering	
	\includegraphics[scale=1]{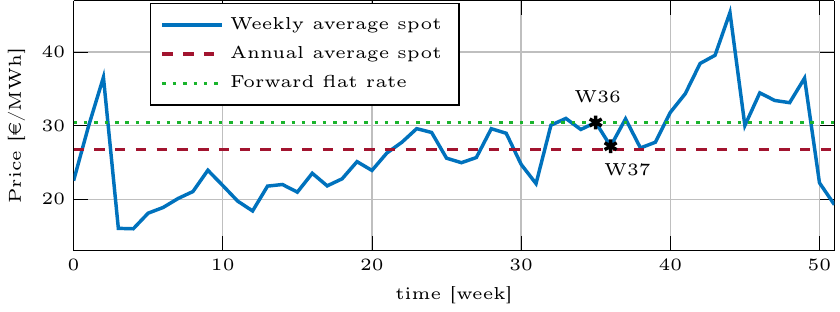}
	\caption{
		Medium weekly price of Nord Pool spot market for DK1 in 2016. Weeks 36 and 37 are highlighted}
	\label{fig:Price1}
\end{figure}%
\txtred{The Spot Market considered in the case study is the Elspot market, in the Nordic Power Exchange, managed by Nord Pool.
Elspot is the world's biggest day-ahead auction-based exchange and it uses Marginal Price System (MPS).} Bidders can submit either price-dependent or price-independent bids for each delivery hour. The latter bid option would allow the consumer BRPs to guarantee their required power via participation in the auction as mere price takers without risking the rejection as a result of offering a price lower than the one set by MPS. Although this renders a precise price forecasting unnecessary for bidding in the spot market as a consumer, a relatively accurate prediction of price profile pattern is needed to shift the power away from price spikes. Most of price prediction techniques offered in literature are capable of offering such level of precision\\ 
\Cref{fig:Price1} presents the weekly average price of DK1 for the year 2016 from Nord Pool \cite{Nordpool} along with the average price and flat rate forward contract for the same year which is 30.4 \euro/MWh for industrial units with annual consumption of 150 MWh and more \cite{Els}.
\begin{figure}[!t]
	\centering	
	\includegraphics[scale=1]{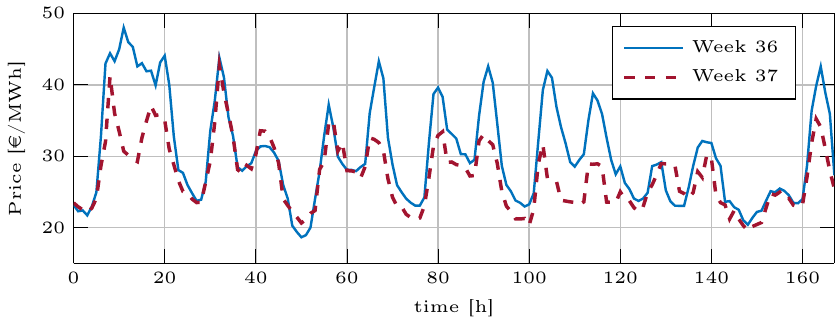}
	\caption{
		DK1 day-ahead hourly price for weeks 36 and 37 of 2016 }
	\label{fig:PriceWeek}
\end{figure}%
We simulated the participation in Elspot market using prices from weeks 36 and 37, whose average values (as shown in \cref{fig:Price1}) are respectively equal to the annual flat rate and average spot price. The hourly price for these time periods is given in \Cref{fig:PriceWeek}. We have also compared the response of the model to a Time-Of-Use tariff (TOU) offered for users with demand greater than 100 MWh per year. This tariff is given in \Cref{tab:TOU}.\\
\begin{table*}[!bht]
	\small
	\centering
	\caption{Time Of Use prices. Tariff by the Danish retailer DONG Energy Eldistribution }
	\label{tab:TOU}
	\begin{tabular}{lll|l|l}  
		\toprule
		&        & \multicolumn{1}{c}{Off-peak}                                                             & \multicolumn{1}{|c}{Shoulder}                      & \multicolumn{1}{|c}{Peak}                          \\ \midrule
		Tariff  [\euro/MWh]                &        & \multicolumn{1}{c}{\cellcolor[HTML]{EFEFEF}27.04}                                        & \multicolumn{1}{|c|}{\cellcolor[HTML]{EFEFEF}33.15} & \multicolumn{1}{c}{\cellcolor[HTML]{EFEFEF}39.39} \\ \midrule
		& Summer &     & 06:00-08:00, 12:00-21:00    & 08:00-12:00     \\ \cmidrule(lr){2-2} \cmidrule(l){4-5} 
		\multirow{-2}{*}{Time slot [hour]} & Winter & \multirow{-2}{*}{\shortstack[l]{21:00-06:00,\\ week-ends, holidays}} & 06:00-08:00, 12:00-17:00, 19:00-21:00  & 08:00-12:00, 17:00-19:00    \\ 
		\bottomrule
	\end{tabular}
\end{table*}%
\txtred{The optimization model was solved with the commercial solver Gurobi 7.0.1 in an Intel Core i7 machine at 2.40 GHz with 16 GB RAM. The average CPU time for the day-ahead scheduling reported in this section was 1185s, subsequently reduced to a total of 351s using the decomposition and warm-start technique proposed in appendix \ref{app:a1}, which are both well within the bidding time window.}
For the sake of comparison, here we define as a measure the Equivalent Flat Rate (EFR), which is the price at which a long-term forward contract would cost exactly the same as participation in the day ahead market. The reason for definition of EFR is to highlight benefits of the proposed DR comparing to a forward contract without revealing sensitive information about the plant under study.
The \Cref{tab:ERF} reports the resulting EFR for the selected weeks.\\
\begin{table*}[!htb]
	\small
	\centering
	\caption{Equivalent Flat Rate for optimal scheduling in day-ahead market. The weekly mean prices, the average of the spot market profiles shown in \cref{fig:PriceWeek}, are compared with the EFR computed each day for the DA optimal solution. \textit{Last week} indicates that the optimal scheduling has been computed using day-ahead price profiles of the previous week. All the entries are expressed in \euro/MWh }
	\label{tab:ERF}
	\begin{tabular}{clc|cccccccc}
		\toprule
		 & \multicolumn{2}{c|}{Day-ahead Price}  &  \multicolumn{8}{c}{Equivalent Flat Rate (EFR)}  \\  \midrule                                                                  
		   & Estimation   &   Weekly Mean   & Mo     & Tu     & We     & Th     & Fr     & Sa     & Su     & Average                            \\ \midrule
		\multicolumn{1}{l}{Week 36}   & \multicolumn{1}{l}{\textit{Exact}    }   & \cellcolor[HTML]{EFEFEF}30.44       & 35.41 &28.63 &27.46 &30.12 &30 &26.11 &25.48 & \cellcolor[HTML]{EFEFEF}29.03 \\ \midrule
		& \textit{Exact    }                     & \cellcolor[HTML]{EFEFEF}27.2      & 29.09 & 28.24 & 26.52 & 26.7 & 24.04 & 25.2 & 23.64 & \cellcolor[HTML]{EFEFEF}26.2 \\
		 \cmidrule(l){2-11} 
		\multirow{-2}{*}{Week 37} & \textit{Last week}  &                                 &30.29 &28.29 &26.75 &26.78 &24.28 &25.35 &23.67& \cellcolor[HTML]{EFEFEF}26.49 \\ \bottomrule
	\end{tabular}
	
\end{table*}
\begin{figure*}[!ht]
	\centering
	\includegraphics[scale=1]{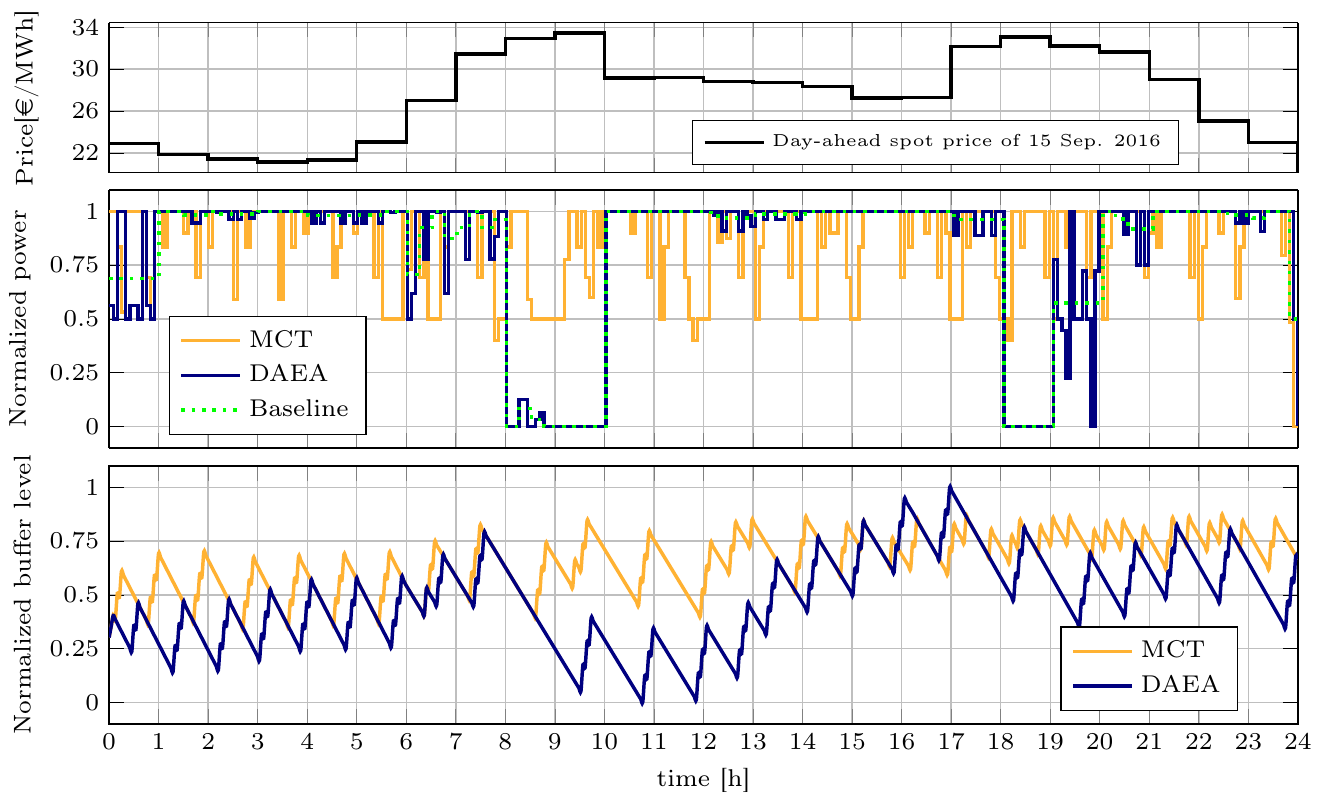}
	\caption{\label{fig:DAEA_P}
		Day-Ahead Energy-Aware scheduling (DAEA, in blue) compared with Minimum Cycle Time (MCT, in ochre) scheduling. Upper panel: DA price profile. Intermediate panel: power profile of four melting IF. Baseline on $\Delta T^{\rm{bl}}$ grid (dashed green). Lower panel: casting buffer level
	}
\end{figure*}%
Besides the fact that a major advantage of being able to purchase the electricity from spot market is its lower average cost (27.2 \euro/MWh comparing to 30.4 \euro/MWh of annual flat rate contract), the resulting EFR for both week 36 and 37 is lower than the average daily price of spot market. According to this outcome, a foundry with a daily consumption of 300 MWh can approximately save up to 400~k\euro/year using the proposed DR approach and participating in Elspot instead of buying electricity from retailers via long-term contracts.\\
Moreover, with respect to a long-term TOU tariff, in which the EAS optimization resulted in an ERF of 27.8~\euro/MWh, participation in day-ahead energy market has proved to be more profitable, even comparing it with average profit of week 36, whose mean spot price is approximately equal to the one of TOU contract. This is mainly due to a single peak hour tariff which hinders the exploitation of process flexibility by the optimizer. While a down-to-zero reduction of consumption level consecutively for a long period may not be possible due to operational and capacity constraints, the same aggregated sum could be obtained if such reduction is broken down to multiple slots well-apart from each other. Accordingly, the optimizer can take advantage of multiple peaks in price and minimize the costly energy consumption.\\
Since this paper does not deal with the issue of price forecasting, for each day in order to verify the effect of error in the prediction of price here we have performed the  optimization also using the price profile of the previous week (W 36), for which the obtained ERF has been reported in the \Cref{tab:ERF}. As it can be noticed, there is just a little discrepancy between the minimal cost obtained using the accurate price and the one minimized by past data, which confirms low sensitivity of the optimizer to precision of forecast price.\\
In order to show the role of flexibility in changing the consumption pattern, the scheduling under DR programme has been compared to one with minimum cycle time objective (MCT), which we assume to be sufficiently similar to the common practice  currently in operation in the foundry under study here. MCT minimizes the cost function \eqref{eq:spot_J}  considering a constant price, therefore due to the time term in energy equation \eqref{eq:loss}, leading to a minimum cycle time. \\
\Cref{fig:DAEA_P} compares the  sum of normalized power profile for four furnaces that deliver to a single buffer. The ochre and blue curves respectively indicating the scheduled consumption for MCT and DAEA. The resultant committed load for the latter case has been shown with the green dashed curve.
From the figure, it is possible to see how consumption pattern has been modified according to price fluctuation, shifting the load  away from time slots with higher energy cost.\\
In the lower panel of  \Cref{fig:DAEA_P}, we have also presented the normalized level of the buffer, using the same colours applied in the power profile panel.
While the MCT solution maintains a level between 25\% and 75\% of the total capacity, DAEA optimization exploits the extremes of buffer capacity, in order to maintain the prescribed outflow even during low-power intervals. Accordingly, a rise in the buffer level before each price peak allows the postponement of cycle completion into cheaper periods, as it can be seen in \Cref{fig:DAEA_P}, where the molten level rises from minimum to maximum in the period between two the peaks.  \\
\begin{figure*}[!ht]
	\centering
	\includegraphics[scale=1]{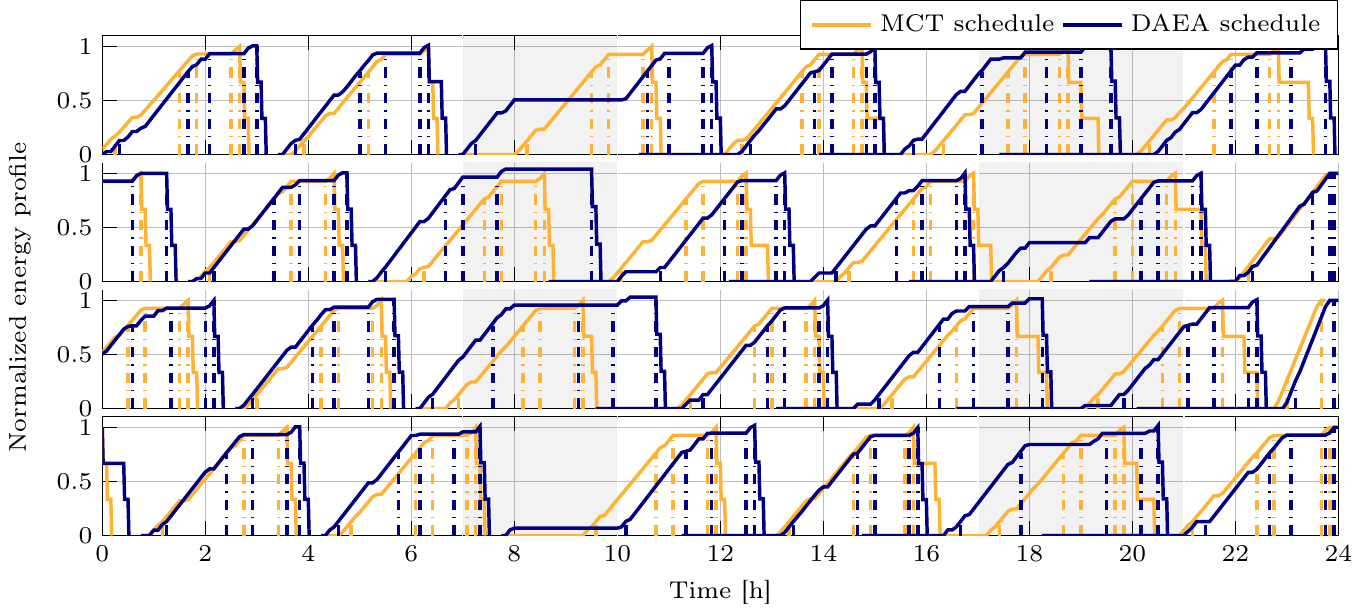}
	\caption{\label{fig:DAEA_E}
		Normalized energy profile of the array of four melting furnaces (DAEA in blue, MCT in ochre). Grey shadowed area are intervals with day-ahead price above $30$~\euro/MWh: DAEA schedule protracts jobs in these region to minimize energy cost}
\end{figure*}%
\Cref{fig:DAEA_E} gives a more detailed explanation of the flexibility provided by the underlying process.
For each furnace, the normalized energy profiles have been depicted for the two scheduling approaches. Vertical dash-dotted lines indicate stage changeovers. 
The grey areas display intervals with spot price over $30$~\euro/MWh. It is evident how DAEA optimization affects the cycles falling within these high-price time slots, while maintaining MCT scheduling for the ones under a relatively flat price.\\
\txtred{It is important to note that, considering the pouring furnace capacity at each casting line, forcing all the first melt cycles to start at the beginning and the last ones to deliver by the end of scheduling horizon would not be optimal. For the same reason, in practice, furnaces delivering to a single line are scheduled to deliver in a sequential manner. Thus, in applying the DR-scheduling model to a running casting line, we assume first cycles to be already started on the previous day, which can be seen in \Cref{fig:DAEA_E}. In addition, by forcing the first and end cycles on each furnace to be complementary to each other, we can impose the same initial condition on the next scheduling to retain the continuity. }
\subsection{Reserve provision}  
The reserve market in Denmark is managed by Energinet.dk and like Elspot is 
divided into western (DK1) and eastern (DK2) bidding zones: regulation varies slightly among the regions. The case study is concerning DK1. It is required to supply the 50\% of the activated reserve within 15 seconds, while the last half must be supplied within 30 seconds. The activated reserve has to be maintained for a maximum of 15 minutes until the secondary and tertiary regulating reserves can take over. 
The bids are eligible to enter the auction if they exceed the minimum of 0.3 MW and received by 15:00 on the day before operation. Bids for more than 5 MW will be rejected altogether, as they can be only accepted in their entirety or not at all. For the auction, a marginal price system is implemented to set the primary reserve price.\\ 
As reported in Energinet.dk database \cite{Energinet}, except for very few hours, the payment rate for primary upward regulation in 2016 remains above 4.4~\euro/\txtred{MW}, while in the tertiary (manual) reserve market only 0.08\% of upward bids have been compensated with more than 0.13~\euro/\txtred{MW}. This confirms the high profitability of primary reserve market with respect to the manual reserve market, particularly considering the much longer sustain time of the latter, which can be as long as 3 hours. Nevertheless, size of the primary reserve market is much smaller, which increases the probability of bid being rejected for non-merit reasons.
Another limiting requirement is the time granularity of the bid. In DK1 the bids should be offered in six equally sized blocks of four hours, while in DK2 the auction for the primary reserve is performed on an hourly basis.\\ 
\begin{figure*}[!h]
	\centering
	\includegraphics[scale=1]{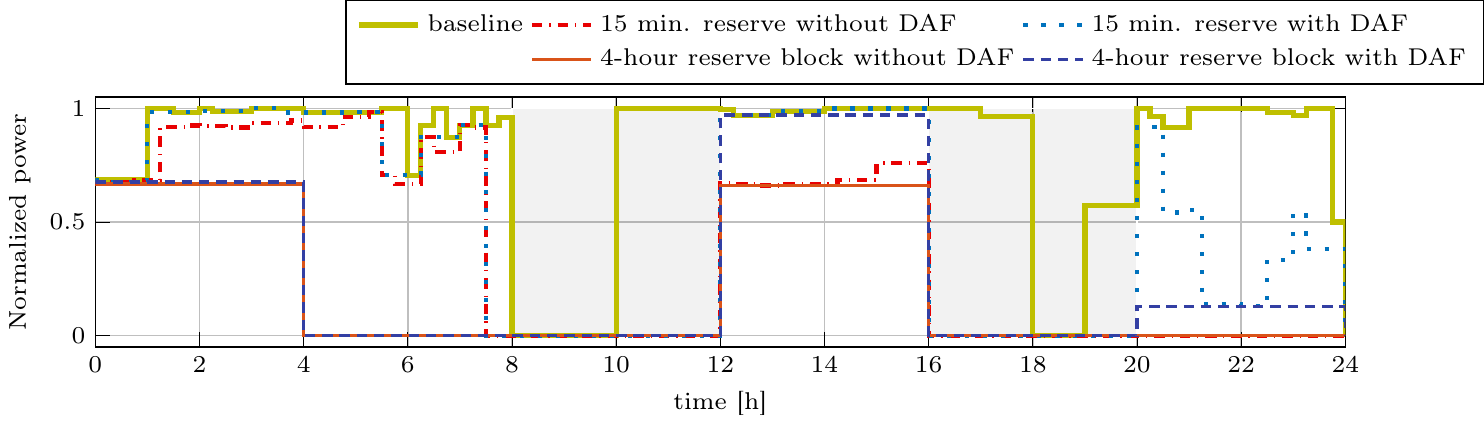}
	\caption{\label{fig:R1}
		Reserve capacity profiles, with (without) day-after flexibility (DAF), are shown in blue dotted (red dash-dotted) line. The 15 minutes profiles are computed with $|\mathcal{Q}|$ independent Reserve Scheduling Optimizations with respect to the Baseline (green) obtained by DAEA scheduling. The 4-hour reserve bids are shown with dashed (solid) lines. }
\end{figure*}%
Unlike the spot market, the behaviour of ancillary and regulating markets cannot be easily predicted. This should not come as a surprise considering the nature of these markets, which is to handle unexpected events causing fluctuations. Many studies have proposed promising probabilistic forecast methods but there is still no model capable of producing informative day-ahead point forecast for efficient balancing market like the Nordic one \cite{klaeboe2015benchmarking}.\\
In the absence of a prediction model here we are going to consider a reasonable worst-case scenario. 
For availability payment, we use a minimum flat rate which according to the past data follows a seasonal pattern. For example, from the beginning of June to the end of September 2016 95\% of the accepted upward bids have been remunerated with a price of 10~\euro/\txtred{MW} while this percentage falls under 25\% for early months of the same year.\\ 
Regarding the imbalance cost, to consider the worst case, we implement a two-price model in spite of the fact that retailer uses a one-price model to determine the price of balancing power for consumption. In addition, we consider a 20\% over/underestimation margin for up/down regulating prices when using the historical data. Even assuming both regulating prices to be equal to the spot one, such a band covers 78\% of down-regulating and  71\% of up-regulating prices of the year 2016. This gives a worst-case probability considering the direction of real-time market to be completely unpredictable. For the case of actual price falling out of this band, we assume the probability of imbalance occurring at the opposite direction of real market to be equal to the probability of their occurrence with the same direction of the market.\\
Such a conservatism needs to be considered in the context of a multi-objective optimization. In this regard, the pricing mechanism is a crucial factor in setting prices for the objective function. 

\begin{figure}[tb]
	\centering
	\includegraphics[scale=1]{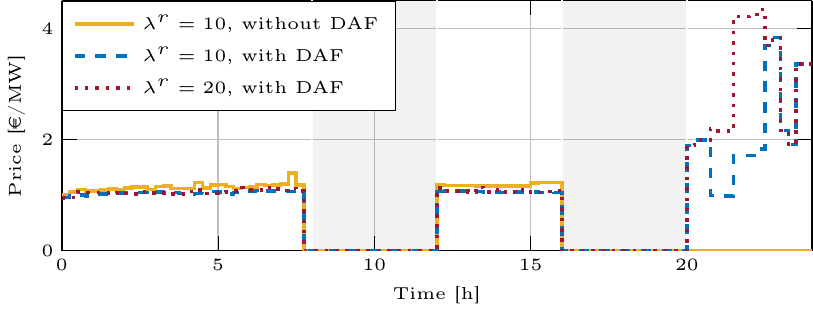}
	\caption{\label{fig:R3}
		Imbalance-induced cost related to 4-hour block bidding intervals for various reserve reward and with and without cycle shift exploiting day-after flexibility (DAF).}
\end{figure}

\begin{figure}[htb]
	\centering	
	\includegraphics[scale=1]{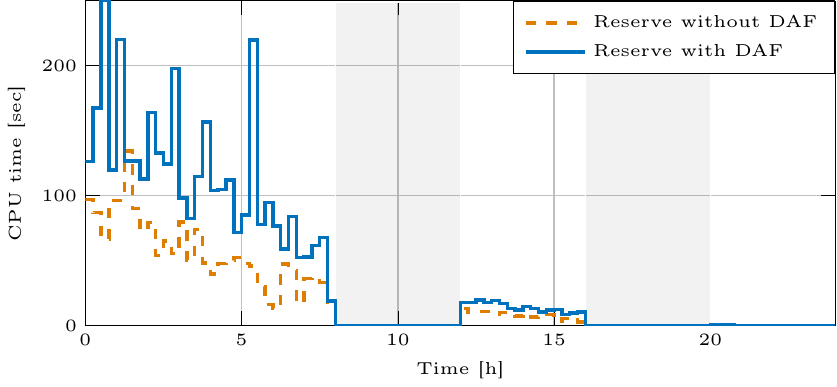}
	\caption{\label{fig:R4_cpu}
		CPU time for each $q \in \mathcal{Q}$ independent Reserve Scheduling Optimizations \eqref{explicitq} with (without) cycle shift for DAF in ochre (blue). On the x axis on behalf  of $q$, the time slot for which the computed reserve $R_{\scriptscriptstyle q,c}$ is related}
\end{figure}

\Cref{fig:R1} shows the offered upward regulation reserve for the same day, whose DAEA results were previously  illustrated in \Cref{fig:DAEA_P}. Here, the grey areas indicate those blocks where bidding is not possible, due to the presence of zero-power instances within them, while the green line shows the spot market committed load profile based on 15-minute settling. 
The calculated reserve for each 4-hour long biding block is shown considering a flat reward rate of 10~\euro/\txtred{MW} for availability,  with above-described 20\% up and down regulating price. Incurred costs have been weighted by the probability factor $ 1/48 $ based on one activation per day with a sustain time of maximum 15 minutes, which renders $ M=2 $ in the optimization problem \eqref{explicitq}. In order to illustrate the advantage of shifting last melting cycles to the next day, we have reported the result of optimization both with and without constraints \eqref{E2next}-\eqref{eq:DeltaE_c}. In the former model, maximum energy price of the day $(D+1)$, $ \bar{\lambda} $, is set to 100 \euro/MWh. This high price has been chosen to reflect the risk associated with the consequent infeasibility of the next day scheduling.\\ 
\begin{figure}[!b]
	\centering	
	\includegraphics[scale=1]{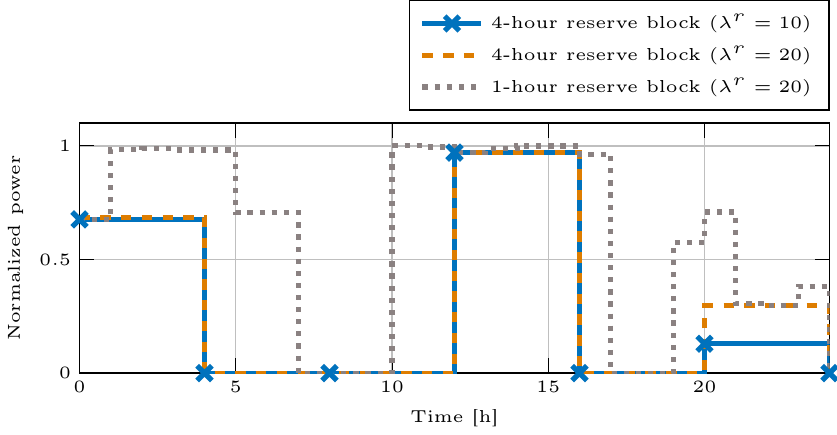}
	\caption{\label{fig:R2}
		Reserve capacity profiles for varying reserve reward and different bidding block hours}
\end{figure}%
\Cref{fig:R2} illustrates the effect of increasing the estimated availability reward for the reserve provision. Moreover, to highlight the limitation of bidding blocks temporal size, we have also included the result for a 1-hour bidding resolution, similar to the regulation of the DK2 market. This result shows a significant increase in the reserve offered under the same condition. However, this otherwise untapped capacity can be exploited using aggregation within the company or via contractual agreement with an external aggregator.\\
Minimum bidding price (i.e. total incurred cost) per reserve unit, for different scenarios, are presented in \Cref{fig:R3}. In all of these cases, the cost is well below the marginal price for reserve availability, which is above 10~\euro/\txtred{MW} during the entire bidding horizon. While this minimum price results in maximal reserve exploitation  at a very low cost for hours prior to 20:00, the twofold increase of the price yields a greater reserve capacity for the last bidding block without inducing a higher maximum cost per unit.\\
Such large profit margin, which is fully attainable within a marginal price auction system, has been achieved in spite of many conservative assumptions we have made regarding the activation cost and duration. For example, in practice,  it normally takes 2-3 minutes for the secondary reserve to take over while we have considered a 30-min interruption interval (for a reserve provision with 15-min resolution). Such a short activation period, resulting in both low net energy imbalance and low probability coefficient in related objective terms, justifies the maximal reserve exploitation using a fictitiously high payment.\\
In the rare event of such conservatism not sufficing, it is easy to see how any subsequent cost could be covered by the generated gross margin.
On the other hand, in a pay-as-bid auction system, it would have been more challenging to ensure a similar benefit margin without a precise estimation of the reserve availability payment rate.\\
Another important requirement that has to be respected for successful market participation is the time for bidding submission. \Cref{fig:R4_cpu} reports the computational time of the optimization performed \txtred{separately for each instant of the reserve grid, with and without day-after flexibility}. Even though the output of two models differ only in the last bidding interval (see  \Cref{fig:R1} and  \Cref{fig:R3}), the inclusion of additional flexibility imposes an extra computational cost, particularly during first intervals. Consequently, we would recommend activating the relaxing constraints \eqref{E2next}-\eqref{eq:DeltaE_c} only during the later bidding intervals to speed up the solution without compromising the bid quality. Furthermore, this minimizes the usage of next day (i.e. $ D+1 $) flexibility when the same capacity can be provided using the energy commitments of the day ahead (i.e. $ D $)  .\\
The overall CPU time for this switching model amounts to 1980 seconds. Although sufficiently small for the conclusion of bid optimization within the time window imposed by TSO, it could still be significantly reduced using parallel computing.\\

\section{Conclusion}
In this paper, we have presented a novel framework for the participation of energy-intensive plants in both energy and reserve markets. 
Due to distinctive features of these markets and different economic priorities of plant owners for entering them, two separate demand response optimization problems have been proposed.
The developed paradigm was tailored to foundry operations as an example of multi-stage and multi-line energy-intensive manufacturing process.\\
First, a mixed integer linear program was formulated to consider the flexibility in power-elastic stages and buffering capacity of storage components such that operational and safety conditions are respected.
This model was used for the optimal scheduling of the day-ahead operation with the goal of minimizing consumption cost considering a variable energy cost.
This enables the participation of plant in the day-ahead energy market whose temporal aspects are reflected in the discretized mixed integer linear model.\\
With the aim of bidding in the capacity market, we further extended this model to account for the incurred cost of resulting imbalance, while maximizing the benefit from offering reserve capacity.
Furthermore, the enhanced model offers additional flexibility by deferring tasks at the end of scheduling horizon to the next day, given the reserve activation probability.
Considering all the requirements of reserve market and calibrating against price uncertainty in the imbalance market, this model produces bidding blocks that ensure the worst-case profit margin. This bidding optimization is performed for each production line separately assuming a maximum activation time of once per day. In order to elude this limitation, we proposed an internal aggregation scheme which optimizes the collective reserve capacity of the plant and increases the possible activation events up to the number of production lines.\\
Using the data from a real-life foundry, we verified both schemes, bidding in the Danish energy and reserve market, considering  price-taking condition.\\
The result of bidding in Nord Pool Elspot market using the proposed model revealed a potential reduction of 9\% in average cost, when compared to the annual forward contract. This outcome showed little dependency on the exact prediction of spot price.\\ Subsequently, the day-ahead baseline was used to evaluate the bidding blocks in the primary reserve market. In spite of high uncertainties associated with the imbalance markets, the proposed model demonstrated to be an effective tool for capitalizing on the day-ahead reserve market, as a secondary revenue stream.\\
Recalling the general framework proposed in \Cref{fig:frame}, in order to guarantee the profitability of such demand-response paradigm, a real-time scheduler should be in place for the execution phase. Development of this module is reserved for the future work where we will focus on the stochastic features of the process and reactive scheduling techniques to handle them. 
Furthermore, in an ongoing research, we are investigating how the proposed Demand Response model can be exploited at the aggregation level for an industrial cluster composed of multiple power-intensive production units. %

\begin{table*}[!t]	
	\begin{framed}			
		\printnomenclature[0.7in]		
	\end{framed}	
\end{table*}

\appendix

\gdef\thesection{\Alph{section}}
\makeatletter
\renewcommand\@seccntformat[1]{Appendix \csname the#1\endcsname.\hspace{0.5em}}
\makeatother

\section{} \label{app:a1}
\subsection{Multiple time grids}
The use of multiple time grid in short-term scheduling, as in \citep{merchan2016discrete}, can be used to reduce the problem size. It results in  a much smaller local time-grid set ${\mathcal{K}_{\scriptscriptstyle {f}, {m},{j}}}= \llbracket k_{\scriptscriptstyle {f}, {m},{j}}^{\scriptscriptstyle \max }-k_{\scriptscriptstyle {f}, {m},{j}}^{\scriptscriptstyle \min }\rrbracket $ 
for each specific task $j_{\scriptscriptstyle {f}, {m}}$, where $\llbracket k \rrbracket=\{1, \dots, k\}$.\\

The complementary unreachable region is defined by the most energetically compact solution for all other stages and jobs. 
Consequently the upper and lower of the local time grid for each stage, i.e. earliest start and latest finish time, can be found as:
\begin{multline}
k_{\scriptscriptstyle f,m,j}^{\scriptscriptstyle \min }=\sum\limits_{\scriptscriptstyle m'\in \llbracket m-1 \rrbracket}{\left[ \sum\limits_{\scriptscriptstyle {j}'\in {{J}_{M}}}{\left\lceil \frac{{{{\hat{E}}}_{\scriptscriptstyle f,{m}',{j}'}}}{p_{\scriptscriptstyle {{j}'}}^{\scriptscriptstyle \max }\delta t} \right\rceil }+\sum\limits_{\scriptscriptstyle j\in {{J}_{T}}}{\left\lceil \frac{{{\Delta }_{\scriptscriptstyle f,{m}',{j}'}}}{\delta t} \right\rceil } \right]}\\
+ \sum\limits_{\scriptscriptstyle {j}'\in {{J}_{M}}\cap  \llbracket j-1\rrbracket}{\left\lceil \frac{{{{\hat{E}}}_{\scriptscriptstyle f,m,{j}'}}}{p_{\scriptscriptstyle j}^{\scriptscriptstyle \max }\delta t} \right\rceil }+\sum\limits_{\scriptscriptstyle {j}'\in {{J}_{T}}\cap \llbracket j-1\rrbracket}{\left\lceil \frac{{{\Delta }_{\scriptscriptstyle f,m,{j}'}}}{\delta t} \right\rceil }
\end{multline}%
\begin{multline}
k_{\scriptscriptstyle f,m,j}^{\scriptscriptstyle \max } \hspace{-1pt} = \hspace{-1pt} K  - 
\hspace{-13pt} \sum\limits_{\scriptscriptstyle m'\in \mathcal{M}_f \setminus \llbracket m \rrbracket}{\left[ \sum\limits_{\scriptscriptstyle {j}'\in {{J}_{M}}}{\left\lceil \frac{{{{\hat{E}}}_{\scriptscriptstyle f,{m}',{j}'}}}{p_{\scriptscriptstyle {{j}'}}^{\scriptscriptstyle \max }\delta t} \right\rceil }+\sum\limits_{\scriptscriptstyle j\in {{J}_{T}}}{\left\lceil \frac{{{\Delta }_{\scriptscriptstyle f,{m}',{j}'}}}{\delta t} \right\rceil } \right]}\\
 - \sum\limits_{\scriptscriptstyle {j}'\in {{J}_{M}}\setminus  \llbracket j-1\rrbracket}{\left\lceil \frac{{{{\hat{E}}}_{\scriptscriptstyle f,m,{j}'}}}{p_{\scriptscriptstyle j}^{\scriptscriptstyle \max }\delta t} \right\rceil }+\sum\limits_{\scriptscriptstyle {j}'\in {{J}_{T}}\setminus \llbracket j-1\rrbracket}{\left\lceil \frac{{{\Delta }_{\scriptscriptstyle f,m,{j}'}}}{\delta t} \right\rceil }
\end{multline}%
where $\lceil x \rceil$ is the ceiling operator.
Consequently, when balancing resources across the multiple lines, as in \eqref{Pmaxl} or \eqref{Pmax}, the summation should be performed over $\mathcal{J}_{\scriptscriptstyle}^{\scriptscriptstyle k}=\left\{ {j_{\scriptscriptstyle f,m}}\in {\mathcal{J}_{\scriptscriptstyle}}|k_{\scriptscriptstyle f,m,j}^{\scriptscriptstyle \min }\le k<k_{\scriptscriptstyle f,m,j}^{\scriptscriptstyle \max } \right\}$  instead of $ \mathcal{J}_{\scriptscriptstyle} $.\\
Moreover, the modeling of the transportation ladle has to be modified by substituting the binary variable in \eqref{eq:3tap} and \eqref{eq:ladle} with:
\begin{equation}
{\tilde{x}}_{\scriptscriptstyle f,m,j}^{\scriptscriptstyle k}= 
\begin{cases}
0 & \textrm{if }  k<k_{\scriptscriptstyle f,m,j}^{\min}, \\
\tilde{x}_{\scriptscriptstyle f,m,j}^{\scriptscriptstyle k} & \textrm{if } k_{\scriptscriptstyle f,m,j}^{\scriptscriptstyle \min }\le k<k_{\scriptscriptstyle f,m,j}^{\scriptscriptstyle \max },\\ 
1 & \textrm{otherwise }
\end{cases}
\end{equation} %

\subsection{Decomposition}

In case a strict time constraint is imposed, e.g. by the market closure or an agreement with BRP, computational time may become an issue.\\
As previously mentioned, the decision about scheduling time step should be seen as a trade-off between computational cost and the solution quality (i.e. optimality). 
However, over tilting this trade-off in favour of computational tractability could result in very conservative or even infeasible solutions.\\
An alternative remedy for speeding up the solution of the optimization problem is a further decomposition at power-pack level, which is possible only for power packs that are entirely associated to a single production line.
\begin{equation*}
\left\{ \exists \bar{\mathcal{L}} \subseteq \mathcal{L} \quad | \quad \bigcup\nolimits_{l \in \bar{\mathcal{L}}} \mathcal{F}_{l}  = \mathcal{F}_{\bar{c}} \right\}
\end{equation*}
To this end, the complication constraint \eqref{eq:buffer} is decoupled by the proportional assignment of buffer volume to furnace $ f \in \mathcal{F}_{l} $.
This may also require the splitting of equation \eqref{buffer_power} and identification of a new $ \gamma $. This possibly sub-optimal solution is then used as a warm-start for the original optimization problem to speed up the overall computational time. In this way, when dealing with very large-scale problems, it is possible to guarantee a feasible (though not optimal) solution, even if optimization has to be interrupted due to time constraint.%
\section*{Acknowledgements}
This work was partially supported by the SYMBIOPTIMA project. SYMBIOPTIMA has received funding from the European Union's Horizon 2020 research and innovation programme under grant agreement No$^{\circ}$  680426
\section*{References}
\bibliography{Biblio} 
\end{document}